\documentclass[preprint,12pt,authoryear]{elsarticle}

\usepackage{amssymb}
\usepackage{adjustbox}
\usepackage{graphics}
\usepackage{amsmath}
\usepackage{subcaption}
\usepackage{graphicx}
\usepackage{booktabs}
\usepackage{xcolor}
\usepackage{comment}
\usepackage[toc,page]{appendix}
\usepackage[utf8]{inputenc}
\usepackage[T1]{fontenc}
\usepackage{graphicx}
\usepackage[hyphens]{url}
\usepackage{hyperref}
\hypersetup{breaklinks=true}
\usepackage{cite}
\usepackage[a4paper, total={6.5in, 9in}]{geometry}
\usepackage{algorithm}
\usepackage{algpseudocode}

\journal{.}

\begin{document}

\begin{frontmatter}

\title{An Improved Metaheuristic Algorithm for On-site Workshop Availability Cost Problem}

\author[a,d]{Niloufar Mirzavand Boroujeni}
\author[b,d]{Nima Moradi\footnote{Corresponding author: Concordia Institue for Information and Systems Engineering, nima.moradi@mail.concordia.ca}}

\affiliation[a]{organization={University of Minnesota},
            addressline={Department of Industrial and Systems Engineering}, 
            city={Minneapolis},
            postcode={MN 55455}, 
            country={USA}}
\affiliation[b]{organization={Concordia University},
            addressline={Concordia Institute for Information and Systems Engineering}, 
            city={Montreal},
            postcode={QC 1455}, 
            country={Canada}}

\affiliation[d]{country={These authors equally contributed to the present work and shared the first authorship.}}

\begin{abstract}
The Multi-mode Resource Availability Cost Problem (MRACP) optimizes resource availability to minimize usage costs and is a recent project scheduling problem variant. The Multi-Mode On-Site Workshop Availability Cost Problem (MOSWACP) extends MRACP by introducing spatial constraints for On-site Workshops (OSWs) at construction sites. MOSWACP aims to determine the optimal availability level, installation, and dismantling times for OSWs while scheduling activities within spatial and resource limitations. A novel Mixed-Integer Linear Programming (MILP) model is developed, and a metaheuristic algorithm, the Electron Radar Search Algorithm (ERSA), is proposed to solve large-scale instances. ERSA, enhanced with problem-specific operators, outperforms CPLEX in large instances and outperforms Simulated Annealing (SA), Genetic Algorithm (GA), and Particle Swarm Optimization (PSO). An actual case study demonstrated significant cost savings using the proposed model. The results and conclusions highlight the effectiveness of the ERSA approach in managing complex project scheduling challenges.

\end{abstract}

\begin{highlights}
\item A new model in project scheduling is presented: Multi-Mode On-Site Workshop Availability Cost Problem (MOSWACP).
\item MOSWACP, formulated as a linear mathematical model, optimizes on-site workshops' lifetime and availability level, beginning time, and execution (operation) mode of activities.
\item An Electron Radar Search Algorithm is developed and enhanced with problem-specific improvement operators.
\item Comparative evaluation with SA, GA, and PSO is conducted alongside implementing the proposed model in a real case study.

\end{highlights}


\end{frontmatter}

\section{Introduction}

A Project Scheduling Problem (PSP) is a combinatorial optimization problem, finding activities' start and finish time under precedence relations constraints to minimize an objective, e.g., timespan and costs \citep{demeulemeester2006project,herroelen2005project}. PSP provides crucial assistance for various projects, like construction projects, to minimize the delays and unused resources caused by inappropriate scheduling. Also, poor project scope determination and design changes cause delays in over 50\% of projects \citep{visual}. Thus, PSP is a tool to reduce these delays and optimize the utilization of resources while minimizing the project's makespan and costs with a direct impact on time, cost, scope, and quality of the project outcome \citep{global}. 

Moreover, resource unavailability, delay of material procurement, space capacity limitation, and project deadline are preliminary sources of complexity for PSP. Variants of PSP are introduced in the literature to involve these sources of complexities. For example, Resource-Constrained PSP (RCPSP) involves resource unavailability, which imposes a limitation on the resource of each activity at different time frames \citep{brucker1999resource,herroelen1998resource,pellerin2020survey,hartmann2022updated}. Time-Constrained PSP (TCPSP) \citep{guldemond2008time}, a.k.a. Resource Investment Problem (RIP) \citep{drexl2001optimization}, is another variant of PSP, which focuses on the project's deadline at minimized resource utilization cost. RIP can also be referred to as the Resource Availability Cost Problem (RACP), which extends RIP to penalize the project's delay or tardiness. RACP minimizes the penalties for project deadline violations instead of minimizing project makespan.

\begin{figure}[t]
    \centering
    \includegraphics[scale=0.35]{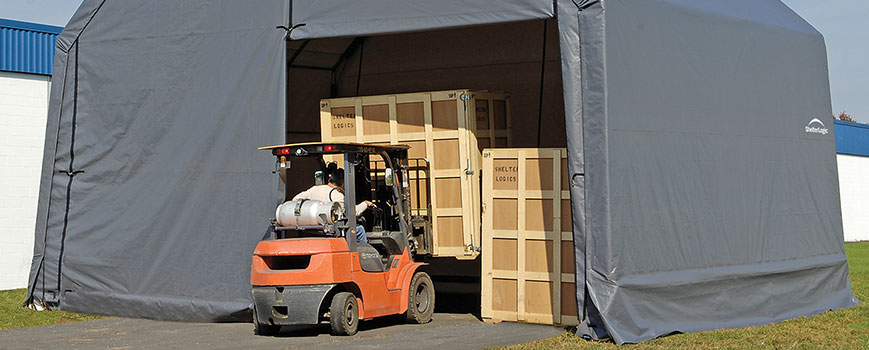}
    \caption{An on-site workshop for storage \citep{Doe:shel:Online}}
    \label{fig:osw}
\end{figure}

Furthermore, the resources in the RACP can be categorized into (non)renewable, doubly constrained, and partially (non)renewable. Renewable resources such as labor, machines, and equipment have limited usage over periods and are replenished for the following periods. In contrast, non-renewable resources, such as money and raw materials, are consumed and not replenished over the planning horizon. Doubly-constrained resources such as periodic budgets are renewable from one period to another but not at each period. Moreover, Partially (non)renewable resources are available during specific periods, and the project activities that require such resources should utilize them only during their availability \citep{bottcher1999project}.
An On-Site Workshop (OSW) is a resource that is partially (non)renewable. involving spatial constraints with a pre-defined lifespan from installation to dismantling. We illustrate an example of OSW in Fig. \ref{fig:osw}, a temporary workshop similar to a garage, protecting project materials vulnerable to heat and rain. OSWs are installed once a project activity requires them and are dismantled upon completion of that activity \citep{Doe:shel:Online}. They effectively reduce the project's costs and makespan by distributing the resources efficiently throughout a project construction or production site considering the project's constraints. 

This paper introduces the Multi-Mode On-Site Workshop Availability Cost Problem (MOSWACP), in which at least a task (activity) can be undertaken in multiple modes to extend RACP constrained to OSWs' spatial capacity and lifetime. The primary motivation for introducing this problem is that the assumption of resources without spatial and lifetime conditions will result in non-optimal (e.g., sub-optimal) solutions. So, by considering the spatial and lifetime-constrained resources, MOSWACP finds more realistic scheduling of activities and resources (e.g., OSWs) than the existing project scheduling-related problems by incorporating real-world aspects of resources and constrained spatial capacity of the project's site. The contributions of the present work are presented as follows:

\begin{itemize}
    \item Addressing a novel construction project problem, multi-mode resource availability cost problem, considering spatial and installation/dismantling times of resources (e.g., OSWs), project site's spatial capacity, and PSP-related constraints. \item Formulating the addressed problem by a novel MILP model with linear constraints. \item Proposing an efficient problem-specific ERSA metaheuristic, comparing with an exact solver and three existing metaheuristics on a new benchmark.
    \item Conducting the comparative study, sensitivity analysis, and a real case study to show the applicability of the proposed model in real situations. 
\end{itemize}

The paper is structured as follows: the related literature is reviewed in Section \ref{sec:background}. Section \ref{sec:model} presents a Mixed-integer Linear Programming (MILP) model for a relatively small-scale MOSWACP. To solve large-scale MOSWACPs, an Electron Radar Search Algorithm (ERSA)-based metaheuristic is proposed in Section \ref{sec:ersa}. Then, the performance of the proposed metaheuristic is compared with a commercial exact solver and existing metaheuristics while conducting a real case study to show the practicability of the proposed model in Section \ref{sec:computation_performance}.
Finally, conclusions, future directions, and model limitations are discussed in Section \ref{sec:conclusions}.

\section{Literature Review}\label{sec:background}

There are several extensions of Resource-Constrained Project Scheduling Problem (RCPSP) in the literature: (a) Multi-Mode RCPSP (MRCPSP) \citep{alcaraz2003solving}, that finds optimal execution mode for each activity; (b) Preemptive MRCPSP \citep{van2010genetic}, which presumes activities could be split or preempted during the project; (c) RCPSP with flexible resource profiles \citep{naber2014mip} in which resource usage of activities can be shifted from one period to another; (d) RCPSP with uncertain resource availabilities \citep{lambrechts2008proactive}; (e) Multi-objective RCPSP \citep{viana2000using} which minimizes the project duration, delay of activities, and resource capacity violation; (f) Stochastic RCPSP \citep{deblaere2011proactive} that involves the uncertain duration of project activities; and (g) MRCPSP with generalized precedence relations \citep{de1999multi} which takes care of a range of time intervals between the start and finish times of the activities.

Furthermore, TCPSP, also known as RIP \citep{drexl2001optimization}, optimizes the availability of resources to minimize resource usage costs while satisfying project deadlines. Several extensions of RIP are introduced in the literature as follows.

\begin{itemize}
    \item RIP with tardiness penalty \citep{shadrokh2007genetic}, which penalizes tardiness of projects.
    \item RIP with discounted cash flows \citep{najafi2006genetic}, which optimizes the net present value of the project cash flows.
    \item RIP with discounted cash flows and generalized precedence relations \citep{najafi2009parameter}.
    \item RIP with time windows \citep{lu2019resource}, in which the periods of resource availability are variant. 
    \item RIP with quantity discount in ordering materials \citep{shahsavar2018integration}.
    \item RIP with discounted cash flows and generalized precedence relations under inflation \citep{shahsavar2010efficient}.
    \item Multi-Mode RIP (MRIP) \citep{hsu2005new,gerhards2018hybrid} in which at least a task can be undertaken in several modes.
    \item Preemptive multi-skilled RIP \citep{javanmard2017preemptive} minimizes the total recruitment cost in the presence of multi-skilled workers.
    \item Uncertain Multi-objective MRIP \citep{subulan2020interval} in which the completion time of activities and availability of renewable resources are uncertain.
\end{itemize}

RACP \citep{zhu2017effective} is another extension of RIP where the tardiness of the project completion time is associated with a penalty. The RACP solvers can be classified into three groups: exact, heuristic, and metaheuristic. Hybrid method with branching scheme \citep{rodrigues2010exact}, Constraint Programming (CP) \citep{kreter2018mixed}, and Modified Minimum Bounding Algorithm (MMBA) \citep{rodrigues2015exact} are among the exact solvers. However, since the exact solvers might be computationally expensive, some heuristic algorithms are developed to solve RACP \citep{peteghem2015heuristic,zhu2017effective,rose2016definition,hsu2005new,su2018constructive,chen2012heuristic}. Also, for RACP with relatively complicated guidelines, metaheuristics are proposed such as Scatter Search (SS) \citep{yamashita2006scatter,meng2016scatter}, SS with a multi-start heuristic \citep{yamashita2007robust}, Artificial Immune System (AIS) \citep{van2013artificial}, path relinking with Genetic Algorithm (GA) \citep{ranjbar2008solving}, pseudo Particle Swarm Optimization (PSO) \citep{qi2014solving}, and Invasive Weed Optimization (IWO) \citep{van2012invasive}.

In addition to the solvers to RACP, several extensions of RACP model are studied in the literature including Multi-Mode RACP (MRACP) \citep{afshar2014multi,qi2014solving,yamashita2009note}, RACP with rental resources \citep{afshar2017project}, RACP with tardiness \citep{su2018constructive}, multi-objective stochastic RACP \citep{arjmand2020evolutionary}, RACP with time-dependant resource cost \citep{afshar2014using}, and RACP with a limited lifetime of resources \citep{afshar2014multi}, from recruitment (installation) to release (dismantling) periods.

In construction or production project scheduling problems, OSWs are partially renewable resources that require physical space, which is limited in construction sites and can be bottlenecks for project completion. We take care of these spatial limitations with efficient OSW installation and dismantling. To this end, we involve OSW in RACP to propose a novel extension to RACP. Solving the project scheduling problems with spatial storage limitations may require expensive computational resources.
Metaheuristic algorithms are developed in the literature to mitigate extensive computations. For instance, \citep{moradi2019simultaneous}, unlike the previous models with single warehouse and unlimited storage capacity, take into account several warehouses with restricted capacity throughout the planning horizon.. They developed a Simulated Annealing algorithm to solve the model. Moreover, \citep{zhang2021project} examined the Project Scheduling and Material Ordering Problem (PSMOP) while accounting for limitation of the storage capacity. The model they developed minimizes the project makespan, stock of materials, ordering costs, and secondary expenses by selecting the optimal activities' timeline, along with determining the timing and amount of material orders.
They also presented an efficient non-dominated sorting GA algorithm to optimize large-sized problems. Moreover, \citep{tian2023integrated} addressed an integrated PSMOP and RCPSP with Limited Storage Space, concentrating on scheduling activities and ordering materials while adhering to storage limitations. They introduced a two-layer heuristic algorithm to solve this model effectively. The problem discussed in the paper, differs from RCPSP-MPS since it involves utilizing esources with partial renewability (or nonrenewability), concurrent activities, and OSWs.

We introduce a problem in project scheduling, MOSWACP, which determines the optimal lifespan and availability threshold of OSWs and the initiation time and execution mode of each activity. Then, we adopt a MILP model for MOSWACP to find optimal solutions for instances of relatively smaller size of MOSWACP and propose an ERSA-based metaheuristic to solve large-scale MOSWACP. Our proposed metaheuristic algorithm enhances ERSA with problem-specific improvement rules, outperforming existing SA, GA, and PSO algorithms. Our research method involves a literature review to understand the problem domain, identifies research gaps, helps adopt a linear mathematical model as a theoretical foundation, and proposes an ERSA metaheuristic to enhance efficiency and performance evaluation by generating different problem instances, including custom-generated benchmarks, and executing tests.

\section{Multi-Mode On-Site Workshop Availability Cost Problem} \label{sec:model}

To formulate MOSWACP, we presume: (1) A single OSW cannot be divided into distinct components; (2) The construction site has restricted space; (3) Each OSW has a unique functionality; (4) The activities have fixed duration times; (5) The activities have multiple execution modes; (6) Each OSW occupies a fixed physical space; (7) Precedence relationships are the finish-to-start (FS) kind; (8) Each OSW is installed at the starting time of the first activity dependent on the OSW. Also, each OSW is dismantled at the finish time of the last activity dependent on the OSW; (9) Tardiness is not permitted; (10) The activities are non-preemptive; and (11) The OSWs are partially (non)renewable resources with pre-determined lifetimes and sizes.

\begin{figure}[t!]
    \centering
    \fbox{\includegraphics[scale=0.6]{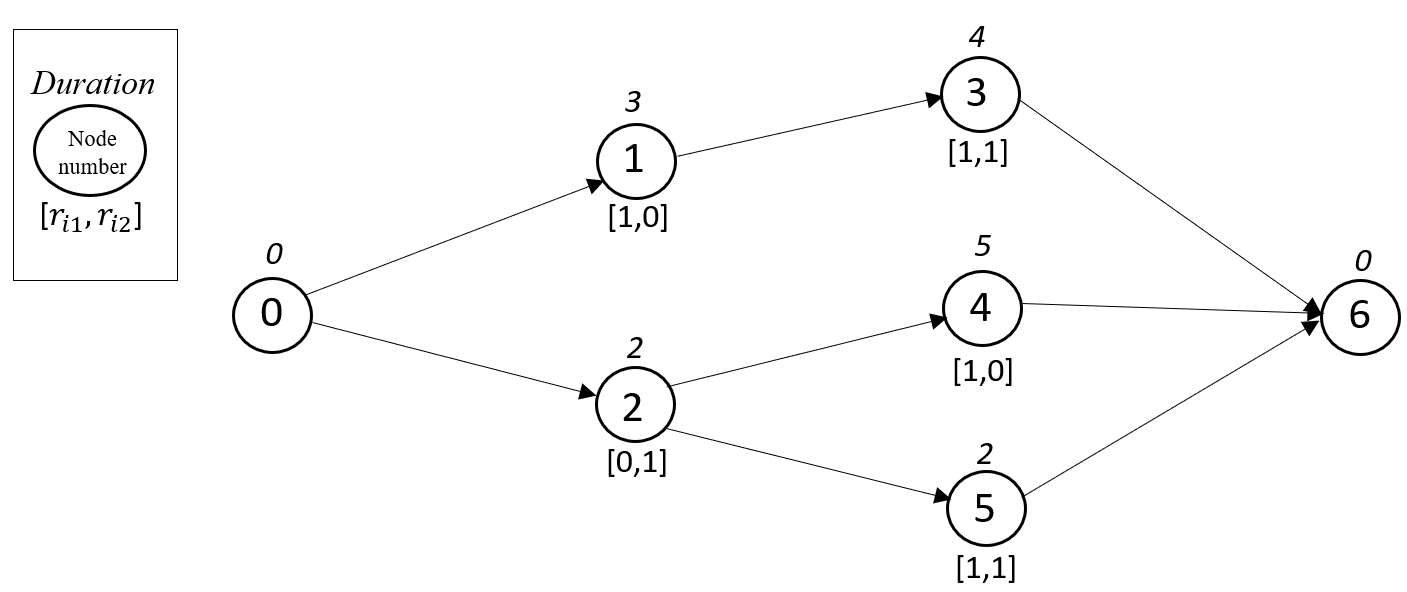}}
    \caption{The example of an AON network}
    \label{fig:graph}
\end{figure}

\begin{figure}[!ht]
\centering
\begin{tabular}{c}
\fbox{\includegraphics[scale=0.5]{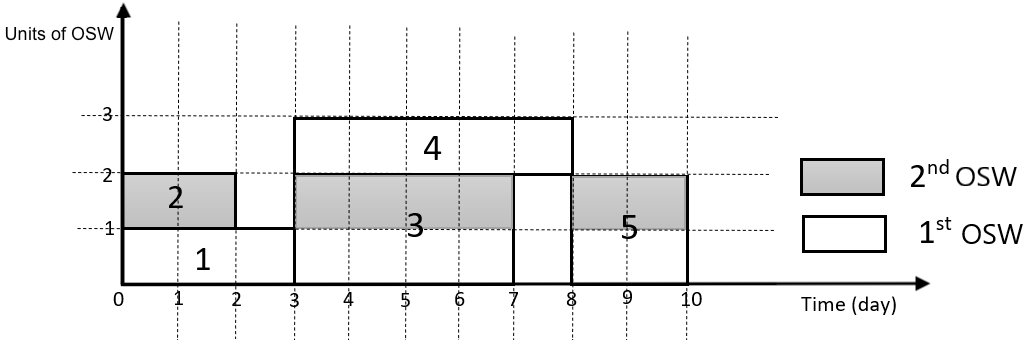}}
\\
\footnotesize (a) A feasible schedule with unlimited physical space
\\ [10pt]
\fbox{\includegraphics[scale=0.42]{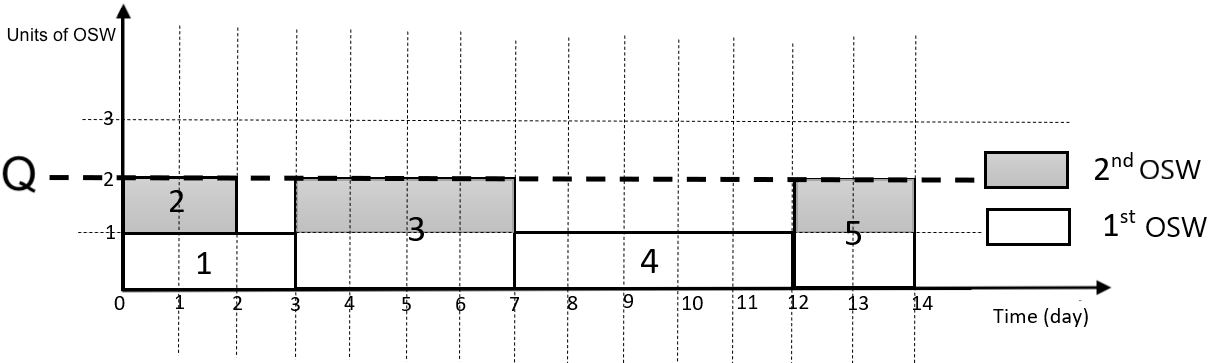}}
\\
\footnotesize (b) A feasible schedule with a limited physical space
\end{tabular}
\caption{The comparison of feasible schedules with (a) limited and (b) unlimited physical space to the AON network in Fig. \ref{fig:graph}}
\label{fig:schedule}
\end{figure}

OSWs significantly influence the construction projects' completion time, cost, and quality. However, due to limited physical space in construction sites, we may be unable to install OSWs at our convenience. Physical space limitations in construction sites may significantly influence the project's timeline. Let's take a construction project with $5$ activities each of which with a single mode of execution, and two OSWs as an example, illustrated by an activity-on-node (AON) network in Fig. \ref{fig:graph}.
We denote the space that the OSW $k$ and activity $i$ require by $r_{ik}$ and add two dummy activities, $0$ and $6$. Also, we specify the duration of each activity right above each activity node in this figure and presume all activities are single-mode types, and precedent relationships are finish-to-start types. A feasible scheduling solution for this project with unlimited and limited physical space is shown in Fig. \ref{fig:schedule}, illustrating the activities' start and finish times and OSW utilization. Even though the project's duration is $10$ days, with limited physical space, imposing a limitation on the physical space delays the project by $40\%$.

We present MOSWACP as a directed graph of the activity-on-node network, $G=(V,E)$, where $V=\{0,1,\ldots,N+1\}$ is the set of activities (nodes) and $E$ is the set of finish-to-start precedence relationship between activities. $V_0$ and $V_{N+1}$ are dummy activities as the initial and terminal nodes, respectively. The activity $j\in V\setminus \{0,N+1\}$ with the execution mode, $i\in m_j$, where $m_j$ is the set of possible modes of execution for activity $j$, with fixed time period, $d_{ji}$. Also, the activity, $j\in J$, with the execution mode, $i\in m_j$, demands a designated physical area, $r_{jik}$, at every OSW, $k\in K$. The precedent activities of activity $j\in J$ are given by $P_j$. The site has a limited area, $Q$, and the cost of utilizing an OSW is $C_k$. The project has a deadline of $T_{max}$, and the planning horizon is $T$. $x_{it}$ is a binary decision variable which equals $1$ if activity $i\in J$ starts at period $t={0,1,\ldots	,T}$ and equals $0$ otherwise. Also, $z_{ij}$ is a binary decision variable that takes $1$ if activity $i\in J$ gets executed on mode $j\in m_i$, and it equals $0$ otherwise. $y_{kt}$ is a binary decision variable and equals one if OSW $k\in K$ is installed at period $t={0,1,\ldots,T}$ and equals $0$ otherwise. $R_k$ denotes the OSW $k$'s level of availability, where $k\in K$.

Using the notations shown in Table \ref{tablevar}, the mathematical programming formulation of MOSWACP, adapted from \citep{moradi2023site}, is presented as Model \eqref{Eq1}-\eqref{Eq10}. The objective function \eqref{Eq1} minimizes the cost of OSWs' usage. Constraint \eqref{Eq2} guarantees the precedence relationships between the activities. Constraint \eqref{Eq3} ensures the occupied area of each OSW does not exceed its level of availability.
Constraint \eqref{Eq4} limits the construction site's physical space capacity.
Constraint \eqref{Eq5} ensures when an OSW-required activity is started, the corresponding OSW is installed and remains active until the activity is finished.
Constraint \eqref{Eq6} imposes a deadline on the project. Constraint \eqref{Eq7} ensures that only one execution mode is selected for each activity. By constraint \eqref{Eq8}, each activity starts only once.
Constraint \eqref{Eq9} sets the initial activity's start time. Finally, the domain of the decision variables is represented by constraint \eqref{Eq10}.

\begin{table}[ht!]
\centering
\footnotesize
\caption{Notations of the model}\label{tablevar}
\begin{tabular}{l p{11cm}}
\hline
\textbf{Notation} & \textbf{Description} \\
\hline
\textbf{Sets} & \\
$J$ & Activities excluding initial and terminal (dummy) activities\\
$m_j$ & Possible execution modes of activity $j\in J$\\
$K$ & Available OSWs\\
$P_j$ & Precedent activities of activity $j\in J$\\
\hline
\textbf{Parameters} &  \\
$N$ & Number of project (non-dummy) activities\\
$d_{ji}$ & Duration of activity $j$ executed on mode $i\in m_j$\\
$r_{jik}$ & The space occupied by activity $j\in J$ on execution mode $i\in m_j$ in OSW $k\in K$\\
$Q$ & The construction site's physical area \\
$C_k$ & Usage cost of OSW $k\in K$\\
$T$ & Planning horizon\\
$T_{max}$ & Project deadline \\
\hline
\textbf{Decision variables} & \\
$x_{it}$ & = 1 if activity $i\in J$ starts at time $t=\{0,1,\ldots	,T\}$; 0 otherwise \\
$z_{ij}$ & = 1 if activity $i\in J$ is executed by the mode $j\in m_i$; 0 otherwise \\
$y_{kt}$ & = 1 if the OSW $k\in K$ is installed at time $t=\{0,1,\ldots	,T\}$; 0 otherwise \\
$R_{k}$ & Availability level of OSW $k\in K$ \\
\hline
\end{tabular}
\end{table}


\begin{eqnarray} \label{Eq1}
&\min. \sum_{K\in k}C_kR_k,&
\\ \notag
&\text{s.t.}&
\\ 
\label{Eq2}
& \sum_{t=0}^{T}tx_{it}\ge \sum_{u\in m_i}z_{ju}d_{ju}+\sum_{t=0}^{T}tx_{jt}, \quad &\forall i\in J, \forall j\in P_i
\\
\label{Eq3}
& \sum_{i\in J}\sum_{j\in m_i}z_{ij}r_{ijk}y_{kt}\le R_{k}, \quad &\forall k\in K, t=\{0,1,\ldots,T\}
\\ 
\label{Eq4}
& \sum_{k\in K}R_{k}y_{kt}\le Q, \quad & t=\{0,1,\ldots,T\}
\\
\label{Eq5}
& tx_{it}+\sum_{j\in m_i}d_{ij}z_{ij}-1 \le \sum_{l=0}^{T}y_{kl}, \quad & \forall i\in J, \forall k\in K,t=\{0,1,\ldots,T\}
\\
\label{Eq6}
& \sum_{t=0}^{T}tx_{n+1,t}\le T_{max}, \quad & t=\{0,1,\ldots	,T\}
\\
\label{Eq7}
& \sum_{j\in m_i}z_{ij} = 1, \quad & \forall i\in J
\\
\label{Eq8}
& \sum_{t=0}^{T}x_{it} = 1, \quad & \forall i\in J
\\
\label{Eq9}
& x_{00} = 1,&
\\
\notag
& x_{it}, z_{ij}, y_{kt}\in \{0,1\}, \quad & \forall i\in J, \forall j\in m_i, \forall k\in K, t=\{0,1,\ldots,T\}\\ \label{Eq10}
\end{eqnarray}

The Model \eqref{Eq1}-\eqref{Eq10} is a nonlinear optimization model with high computational complexity.
To solve the problem with commercial solvers, such as Gurobi or CPLEX \citep{cplex}, we linearize the model by linearizing constraints \eqref{Eq3} and \eqref{Eq4}.

\subsection{Linearization of constraint \eqref{Eq3}}

We define a new binary decision variable, $y'_{ijkt}$. It equals $1$ if activity $i\in J$ with the execution mode $j\in m_i$ occupies the OSW $k\in K$ at the period $t=\{0,1,\ldots,T\}$, otherwise it equals $0$.
Accordingly, we replace constraint \eqref{Eq3} with constraints \eqref{Eq11}-\eqref{Eq13}.

\begin{eqnarray}
    \label{Eq11}
    \sum_{i\in J}\sum_{j\in m_i}r_{ijk}y'_{ijkt}\le R_k,& \quad \forall k\in K, t=\{0,1,\ldots,T\}
    \\
    \label{Eq12}
    y'_{ijkt}+1 \ge z_{ij}+y_{kt},& \quad \forall i\in J, \forall j\in m_i, \forall k\in K, t=\{0,1,\ldots,T\}
    \\
    \label{Eq13}
    y'_{ijkt}\le \alpha (z_{ij}+y_{kt}),& \quad \forall i\in J, \forall j\in m_i, \forall k\in K, t=\{0,1,\ldots,T\}
\end{eqnarray}

Where $\alpha$ is an arbitrary constant parameter such that $0< \alpha <1$.

\subsection{Linearization of constraint \eqref{Eq4}}

We defined a new decision variable, $R'_{kt}$, as the availability level of OSW $k\in K$ at period $t=\{0,1,\ldots	,T\}$. Then, we replace the constraint \eqref{Eq4} with constraints \eqref{Eq14}-\eqref{Eq17}.

\begin{eqnarray}
    \label{Eq14}
    \sum_{k\in K}R'_{kt}\le Q,& \quad t=\{0,1,\ldots,T\}
    \\
    \label{Eq15}
    0\le R'_{kt}\le R_k,& \quad \forall k\in K, t=\{0,1,\ldots,T\}
    \\
    \label{Eq16}
    R'_{kt}\le My_{kt},& \quad \forall k\in K, t=\{0,1,\ldots,T\}
    \\
    \label{Eq17}
    R'_{kt}\ge R_k-M(1-y_{kt}),& \quad \forall k\in K, t=\{0,1,\ldots,T\}
\end{eqnarray}

Where $M$ is a big constant parameter, such that $0\le R_k\le M, \forall k\in K$.

Finally, we reformulate the adopted nonlinear MIP \eqref{Eq1}-\eqref{Eq10} to a linear MIP \eqref{Eq1}, \eqref{Eq2}, \eqref{Eq5}-\eqref{Eq17}, \eqref{Eqx3}.

\begin{eqnarray} \notag
&\min.& \text{Objective function } \eqref{Eq1},
\\ \notag
&\text{s.t.}&
\\ \notag
&&\text{Constraints } \eqref{Eq2}, \eqref{Eq5}-\eqref{Eq17},
\\ \notag
&& y'_{ijkt}, R'_{kt} \in \{0,1\}, \quad \forall i\in J, \forall j\in m_i, \forall k\in K, t=\{0,1,\ldots,T\}.\\ \label{Eqx3}
\end{eqnarray}

The adopted MIP model for MOSWACP has $T\bar{K}m'$ decision variables and $n'+T\bar{K}m'$ constraints.
$\bar{K}$, $m'$, and $n'$ are the total number of available OSWs, the total number of execution modes, and the total number of precedent activities, respectively.
MOSWACP extends the Multi-Mode Resource Investment Problem with tardiness (MRIPT). Since MRIPT is NP-hard \citep{gerhards2020multi}, and MOSWACP can be reduced to MRIPT, the MOSWACP is NP-hard, too. The existing commercial solvers cannot solve a large-scale MOSWACP within a reasonable time. Thus, we developed the metaheuristic algorithm, Electron Radar Search Algorithm (ERSA), to generate solutions for large-scale MOSWACPs.

\section{Electron Radar Search Algorithm for MOSWACP} \label{sec:ersa}

ERSA was first introduced by \citep{rahmanzadeh2020electron}, which mimics the natural behavior of electric flow when the electrons are in a gas, liquid, or pooplid environment. When the voltage between the anode and cathode increases, the electrons emit from a lower to a higher potential. In this situation, the electrons emit through a way by which they encounter the least resistance in the environment. To evaluate the surrounding environment and find the path with less resistance, the electrons use a radar mechanism to search the environment for an efficient path. ERSA finds the best solution in $95\%$ of the benchmark functions \citep{rahmanzadeh2020electron}, so we re-design it to comply with our proposed MOSWACP and solve it. We present the pseudo-code of ERSA in Alg. \ref{fig:1}. In this pseudo-code, $\beta$ controls the balance between exploration and exploitation, $E^0_n$ is the number of electrons for each initial solution, $CV$ is the critical value, $r$ is the searching radius for each solution, and $M$ is the number of random points chosen for the local search at each iteration, $t$.

\begin{algorithm}
\caption{The pseudo-code of ERSA \citep{rahmanzadeh2020electron}}\label{fig:1}
\begin{algorithmic}[1]
\State \textbf{Start}
\State \textbf{End}
\State \textbf{Inputs:} Initial population of streamers ($N$), parameters $\beta$, $E^0_n$, $CV$, $r$, $M$, and $t = 0$
\While {The number of existing streamers is greater than zero}
    \For {each streamer}
        \If {Eliminating condition is satisfied}
            \State Eliminate streamer and go to the next streamer
        \ElsIf {Forking condition is satisfied}
            \State Create a new streamer randomly and update the population of streamers
        \EndIf
        \State Update streamer position and save the best objective function
    \EndFor
\EndWhile
\State Display the best objective function
\end{algorithmic}
\end{algorithm}

\subsection{Solution representation}

We built a $2n+3\bar{K}$-dimensional vector to represent the solution space for the MOSWACP problem as

\begin{center}
     $\{S_1,\ldots	,S_n,SW_1,\ldots	,SW_{\bar{K}},FW_1,\ldots	,FW_{\bar{K}},M_1,\ldots	,M_n,R_1,\ldots	,R_{\bar{K}}\}$.
\end{center}

where $n$ denotes the total number of activities, and $\bar{K}$ represents total OSWs. The solution representation in the proposed ERSA is a ($2n+3\bar{K}$)-dimensional vector, in which $\{S_1,\ldots,S_n\}$, $\{SW_1,\ldots	,SW_{\bar{K}}\}$, $\{FW_1,\ldots,FW_{\bar{K}}\}$, $\{M_1,\ldots,M_n\}$, and $\{R_1,\ldots,R_{\bar{K}}\}$ are the activities' start time, the installation time of the OSWs, the dismantling duration of OSWs, the execution mode selected for activities, and the level of OSWs' availability, in corresponding order.

\subsection{Initial population generation}

The initial population set must constitute a feasible solution for MOSWACP. First, we find $R_k, k\in K$, through the allocation of a randomly generated integer that exceeds $\max_{i\in J, j\in m_i}\{r_{ijk}\}$. Second, we find $M_i (i\in J)$ by selecting a random integer from set $m_i$, which is the execution mode of the activity $i$. Third, we use Alg. \ref{fig:2} to determine the initiation time of the activities. Finally, $SW_k$ and $FW_k$, $k\in K$ are determined based on the activities that require OSW $k\in K$ to begin, as Eq. \eqref{Eq.SW}.

\begin{eqnarray}\label{Eq.SW}
    SW_k=\min_{i: \exists j\in m_i, r_{ijk}\neq 0}\{S_i\}, \quad FW_k=\max_{i: \exists j\in m_i, r_{ijk}\neq 0}\{S_i+d_{ij}\}  
\end{eqnarray}

\begin{algorithm}
\caption{Feasible schedule generator for project activities}\label{fig:2}
\begin{algorithmic}[1]
\State \textbf{Start}
\State Set $S_0 = 0$
\While {The scheduling is not finished}
    \State Choose the next activity in the set $\{S_1, \dots, S_n\}$
\For {$i = 1$ to $n$}
    \State Set $S_i = \max_{u \in P_i} \{S_u + d_{uj}\}$ for $j=M_u$ \Comment{$P_i$ is the set of precedent activities of activity $i$, $M_u$ has already been determined in the chromosome}
    \State $R_k$
    \If {$R_k \leq Q$}
        \State $R_k \gets R_k + r_{ijk}$ for $j = M_i$ \Comment{Updating $R_k$}
    \EndIf
\EndFor
\EndWhile
\State \textbf{Return} $S_1, \dots, S_n$
\State \textbf{End}
\end{algorithmic}
\end{algorithm}

\subsection{Improvement-based local search}

ERSA employs an improvement-based local search to create the neighboring solution by practical crossover and mutation operators.
Then, it is improved by applying three efficient improvement rules represented in Table \ref{table1}. The inputs and outputs of the crossover operator are $2$ solutions, each of which known as parent and offspring solutions. We presume $Sol_1$ and $Sol_2$ are parent solutions as follows.

\begin{center}
     $Sol_1:\{S^1_1,\ldots	,S^1_n,SW^1_1,\ldots	,SW^1_{\bar{K}},FW^1_1,\ldots	,FW^1_{\bar{K}},M^1_1,\ldots	,M^1_n,R^1_1,\ldots	,R^1_{\bar{K}}\}$.
\end{center}

\begin{center}
     $Sol_2:\{S^2_1,\ldots	,S^2_n,SW^2_1,\ldots	,SW^2_{\bar{K}},FW^2_1,\ldots	,FW^2_{\bar{K}},M^2_1,\ldots	,M^2_n,R^2_1,\ldots	,R^2_{\bar{K}}\}$.
\end{center}


We determine a new point based on the crossover point of a single-point crossover, $C_p$.
The crossover then takes place at this point, producing the offspring solutions, $Sol_3$ and $Sol_4$, as described below.

\begin{center}
     $Sol_1:\{S^1_1,\ldots	,S^1_n,SW^1_1,\ldots	,SW^1_{\bar{K}}\otimes ,FW^1_1,\ldots	,FW^1_{\bar{K}},M^1_1,\ldots	,M^1_n,R^1_1,\ldots	,R^1_{\bar{K}}\}$
\end{center}

\begin{center}
     $Sol_2:\{S^2_1,\ldots	,S^2_n,SW^2_1,\ldots	,SW^2_{\bar{K}}\otimes,FW^2_1,\ldots	,FW^2_{\bar{K}},M^2_1,\ldots	,M^2_n,R^2_1,\ldots	,R^2_{\bar{K}}\}$
\end{center}

\begin{center}
    $\downarrow$
\end{center}

\begin{center}
     $Sol_3:\{S^1_1,\ldots	,S^1_n,SW^1_1,\ldots	,SW^1_{\bar{K}}\otimes,FW^2_1,\ldots	,FW^2_{\bar{K}},M^2_1,\ldots	,M^2_n,R^2_1,\ldots	,R^2_{\bar{K}}\}$
\end{center}

\begin{center}
     $Sol_4:\{S^2_1,\ldots	,S^2_n,SW^2_1,\ldots	,SW^2_{\bar{K}}\otimes,FW^1_1,\ldots	,FW^1_{\bar{K}},M^1_1,\ldots	,M^1_n,R^1_1,\ldots	,R^1_{\bar{K}}\}$
\end{center}

With the symbol $\otimes$ and $C_p=n+\bar{K}$, crossover point is positioned between $n+\bar{K}$ and $n+\bar{K}+1$.

Since the offspring solutions might be infeasible, we designed a crossover operator to transform an infeasible solution to a feasible one. If the crossover modify the beginning time of the activities, we call Alg.~\ref{fig:2} to repair the infeasible schedule. Additionally, if the crossover alters the OSW schedule, we transform the infeasible schedule into a feasible one by Eq. \eqref{Eq.SW}.
Moreover, in case the crossover shifts the execution mode of an activity, we update the schedule of the activities and OSW according to the new execution modes by Alg.~\ref{fig:2}. Finally, in case the crossover modifies the OSWs' availability level, we update the schedule of OSWs and activities based off of the new availability levels using Alg.~\ref{fig:2}.

The mutation operator takes one solution, and creates a mutated solution. The number of bits (elements) mutated in each solution is determined by $r_m$. Let's take the solution, $Sol_1$, to mutate as follows.

\begin{center}
     $Sol_1:\{S^1_1,\ldots	,S^1_n,SW^1_1,\ldots	,SW^1_{\bar{K}},FW^1_1,\ldots	,FW^1_{\bar{K}},M^1_1,\ldots	,M^1_n,R^1_1,\ldots	,R^1_{\bar{K}}\}$
\end{center}

We mutate the elements $n$ and $n+\bar{K}$, such that $r_m=2$. The new solution becomes as follows.

\begin{center}
     $\{S^1_1,\ldots	,\mathbf{S'^1_n},SW^1_1,\ldots	,\mathbf{SW'^1_{\bar{K}}},FW^1_1,\ldots	,FW^1_{\bar{K}},M^1_1,\ldots	,M^1_n,R^1_1,\ldots	,R^1_{\bar{K}}\}$
\end{center}
Where $S'^1_n$ and $SW'^1_{\bar{K}}$ are the new elements.

The solution may become infeasible by mutation. In that case, we use Alg. \ref{fig:2} to convert infeasible solutions into feasible solutions. If the mutation changes the activity schedule, we convert the activity's start times to satisfy precedence relationships and the space capacity constraints. Also, if the mutation operator changes the installation or dismantling times of the OSW, we modify the corresponding installation or dismantling times to maintain the feasibility of the activities schedule. Likewise, if the mutation operator modifies the availability level of an OSW, we modify the corresponding OSW's availability level to maintain the feasibility of the schedule and space capacity.
Finally, if the mutation operator shifts the execution modes, we switch that activity to support the feasibility of activities' schedules and space capacity.
As the fitness function involves the cost of OSWs utilization, we calculate total usage cost by $\sum_{k\in K}R_kC_k$.

The proposed ERSA in this study incorporates three MOSWACP-specific improvement operators, as shown earlier in Table \ref{table1}.

\begin{table}[!ht]
\footnotesize
\caption{Improvement rules in the proposed ERSA}
\label{table1}
\centering
\begin{tabular}{p{6.5cm} p{6.5cm}}
\hline
Improvement rule & Description \\
\hline
Switching the execution mode, $IR_1$ & Switch the execution mode of activities while maintaining the feasibility of solutions as well as improving the value of the objective function \\
\hline
Rescheduling activities, $IR_2$ & Reschedule the activity's start time based on its float time and improve the values of the objective function at the same time \\
\hline
Regularizing OSWs' level of availability, $IR_3$ & Regularize the availability of an OSW without violating the physical space limitation an OSW takes \\
\hline 
\end{tabular}
\end{table}

The first improvement operator, denoted as $IR_1$, focuses on altering the activities' mode of execution to achieve two primary goals: (i) Reduce the project's makespan by parallelizing activities, and (ii) Optimize the utilization of OSWs and avoid the physical space violation.

\begin{figure}[ht!]
    \centering
        \fbox{\includegraphics[scale=0.8]{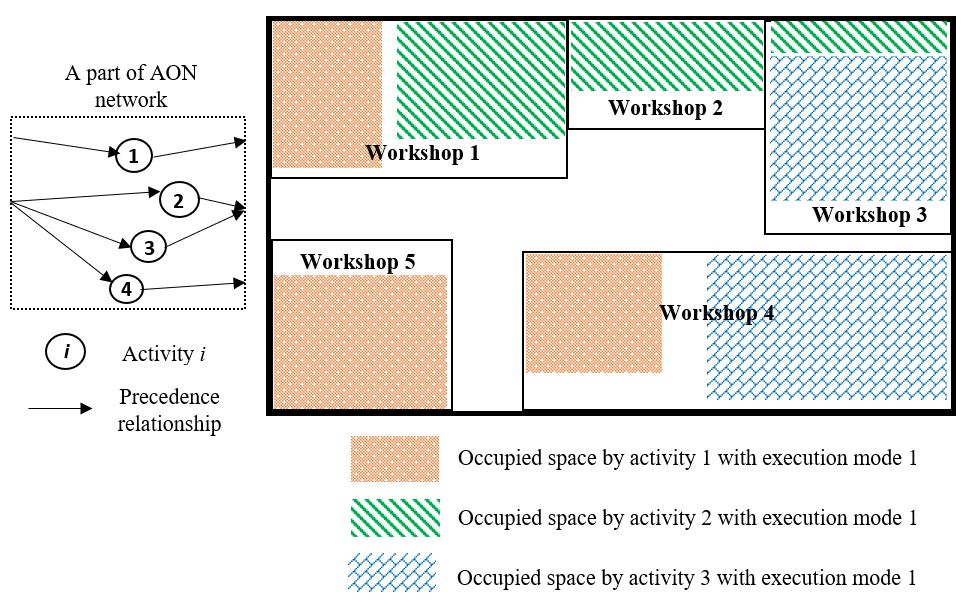}}
    \caption{A construction site with $5$ OSWs and the corresponding project AON network, where activity $1$ is executed on mode $1$}
    \label{fig7}
\end{figure}

\begin{figure}[ht!]
    \centering
    \fbox{\includegraphics[scale=0.8]{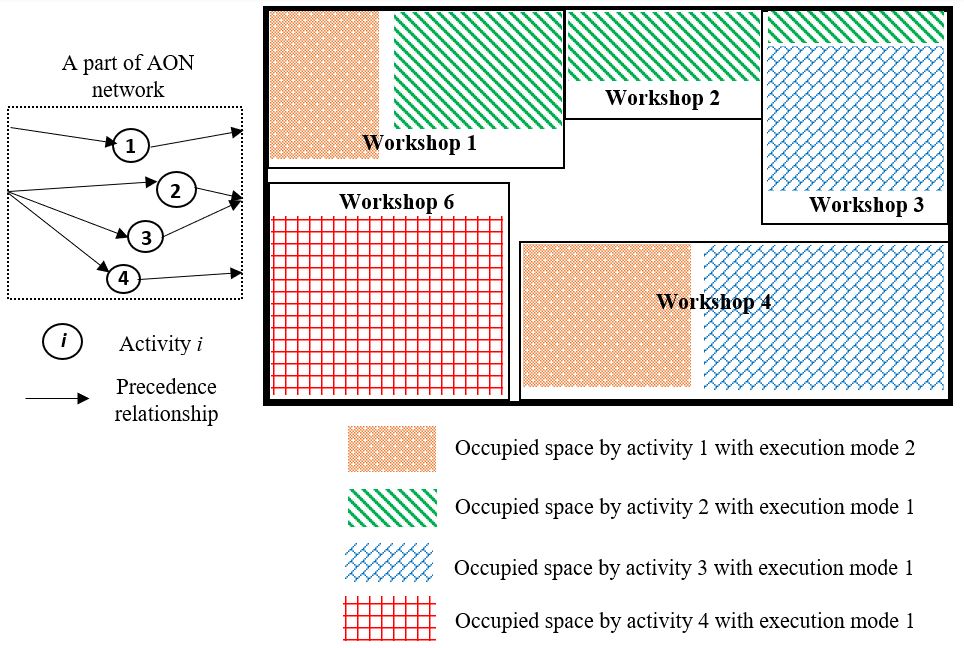}}
    \caption{A construction site with $5$ OSWs and the corresponding project AON network, where the execution mode of activity $1$ is switched from mode $1$ to mode $2$ by $IR_1$}
    \label{fig8}
\end{figure}

Let's take an example of a construction site, shown in Fig. \ref{fig7}, with $5$ OSWs, in which the activity $1$ utilizes the OSWs $\{1,4,5\}$, activity $2$ utilizes $\{1,2,3\}$, and activity $3$ utilizes $\{3,4\}$. The improvement operator, $IR_1$, switches the activity's execution mode $1$ from mode $1$ to mode $2$. In this case, the second execution mode requires activity $1$ to take relatively more space in workshop $4$, while not demanding a space in OSW $5$, as illustrated in Fig. \ref{fig8}. Consequently, we can execute activity $4$ concurrent with activities $1$, $2$, and $3$ due to the precedence relationships outlined in Fig. \ref{fig7}. Activity $4$ occupies workshop $6$ only so that our adjustment does not violate the physical space capacity. In this example, we can reduce the project duration by effectively parallelizing $4$ activities without violating the constraints.

\begin{figure}[ht!]
    \centering
    \fbox{\includegraphics[scale=0.8]{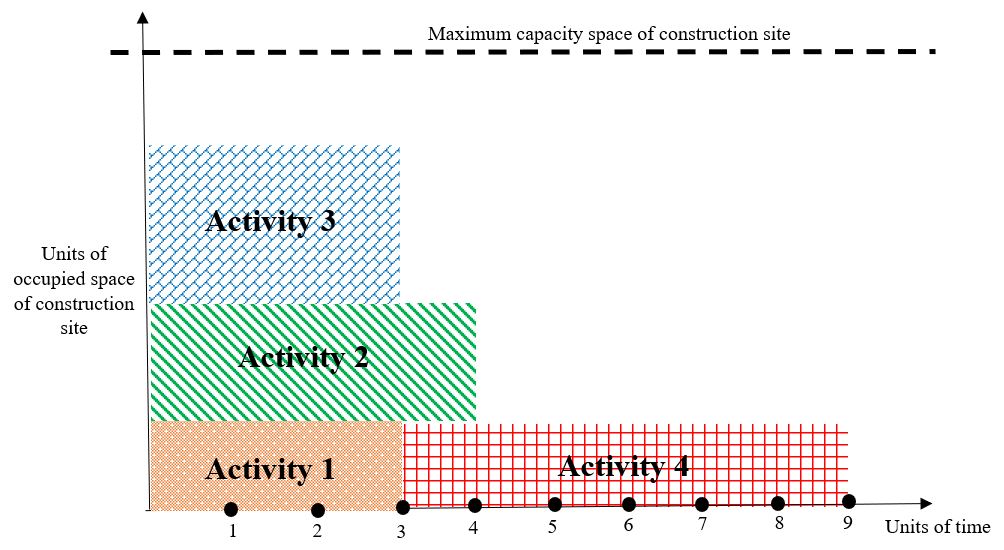}}
    \caption{Activities schedule when activity 1 starts at time zero}
    \label{fig9}
\end{figure}

\begin{figure}[ht!]
    \centering
    \fbox{\includegraphics[scale=0.8]{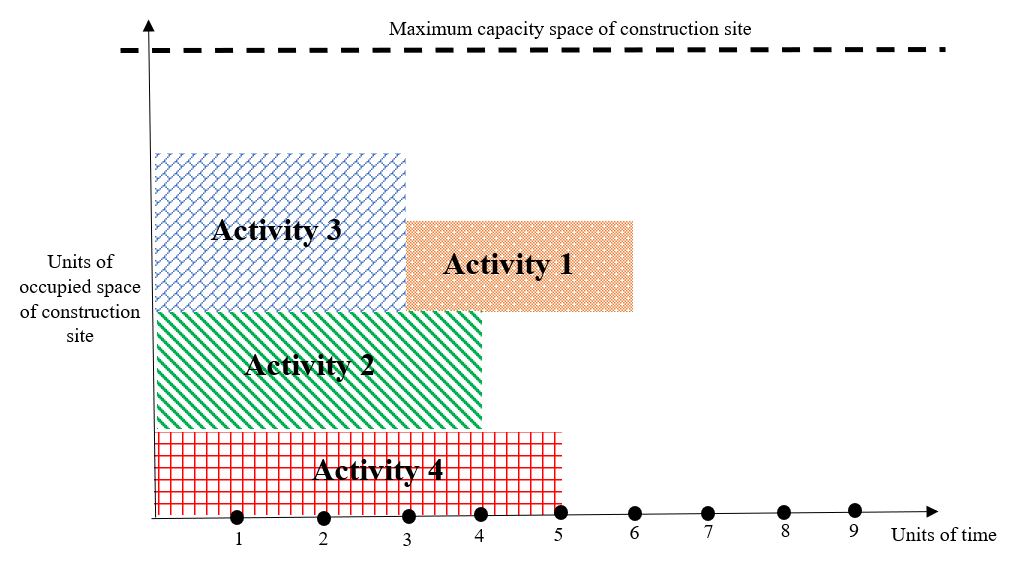}}
    \caption{Activities schedule when activity one is rescheduled to start at time three by $IR_2$}
    \label{fig10}
\end{figure}

The second improvement operator, denoted as $IR_2$, rearranges the scheduling of activities based on their float times and reduces the project's makespan, consistent with precedence relationships and the area's constrained capacity. To illustrate the improvement operator $IR_2$, let's take the scenario depicted in Fig. \ref{fig7} as an example, where the duration of activities $1-4$ are $3$, $4$, $3$, and $6$ periods, respectively. Due to the limited area at the construction site, activity $4$ cannot commence until, say, one of the activities $1-3$ is complete. A feasible schedule for this project is shown in Fig. \ref{fig9}. In this example, activity $1$ is not a direct successor of $2$ or $3$. In such cases, the improvement operator, $IR_2$, shifts the activity's start time $1$ from $0$ to $3$; see Fig. \ref{fig10}. This rule reduces the overall project duration by $33\%$; see Fig. \ref{fig9}, and  Fig. \ref{fig10}.

\begin{figure}[ht!]
    \centering
  \fbox{\includegraphics[scale=0.8]{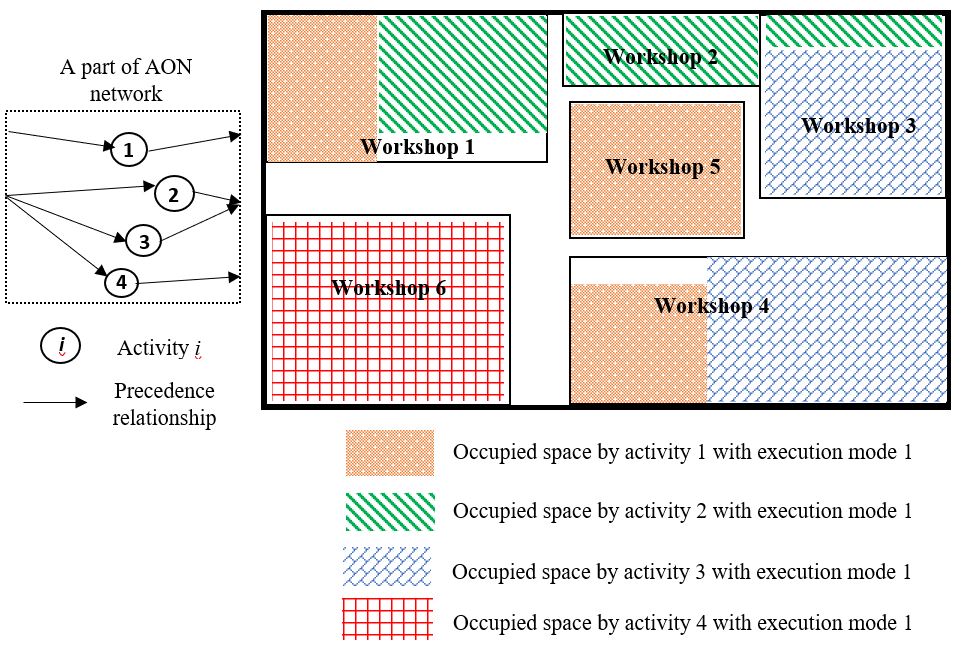}}
    \caption{The construction site, occupied with 6 OSWs while the size of each OSW is regularized by $IR_3$}
    \label{fig11}
\end{figure}

The improvement operator, $IR_3$, optimizes the size or availability level of OSWs to incorporate additional workshops on the site. Let's take the construction site depicted in Fig. \ref{fig7} as an example. The improvement operator, $IR_3$, reduces OSWs' availability level and maintains the demanded physical space of activities in them. By resizing the workshop(s), we add space to the site and install an additional OSW, as Workshop $6$ in Fig. \ref{fig11}. This adjustment facilitates the concurrent execution of more activities and ultimately improves the project schedules and allocation of resources efficiently.

\section{Computational Results} \label{sec:computation_performance}

In this section, we solve the MOSWACP with MIP solver, ERSA, and well-known metaheuristics, including SA \citep{moradi2022efficient}, GA \citep{moradi2023site}, and PSO \citep{kayvanfar2023hybrid}. Also, we compare the performance of the solvers. We solve the MIP model using the CPLEX Application Programming Interface (API) in Python and code the ERSA in C\verb!++! on a 1.60 GHz Core i5 processor with 16GB of RAM. The parameters of MOSWACP are close to the well-known MRCPSP \citep{alcaraz2003solving}. To generate the instances for MOSWACP, we add the utilization cost of each OSW, $C_k$, and its corresponding demanding area utilized with execution mode $i$, $r_{jik}$ of activity $j$, to the MMRCPSP instances.

\subsection{Instance generation for MOSWACP}

In this section, we explain how we involve novel parameters used in the literature to generate diverse instances of MOSWACP. The parameters for MOSWACP closely resemble those of the well-known MMRCPSP \citep{drexl1993nonpreemptive}. To create MOSWACP instances, we involve $C_k$ and $r_{jik}$ in MMRCPSP library\footnote{http://www.om-db.wi.tum.de/psplib/getdata\_mm.html}, which represent the usage cost of the OSW, $k$, and the space occupied by activity, $j$, with execution mode, $i$, respectively.
We selected the instances, \textit{C15, C21, J10, J12, J14, J16, J18, J20, M2}, and \textit{R3}, from this library. We evaluate how the problem responds to the parameter variations using $16$ different instances of \textit{J10}. Subsequently, we incrementally range the parameter $r_{jik}$ from $0$ to $13$, with other parameters fixed. Similarly, incrementally, we set the parameter $C_k$ to several different values from $\{10,20,30,50,100\}$ with other parameters fixed. Finally, we generate MOSWACP instances of various sizes, 34, 84, 536, 101, 117, 45, 93, 64, 89, and 113, for the datasets C15, C21, J10, J12, J14, J16, J18, J20, M2, and R3.

\subsection{Parameter tuning}\label{param_tun}

We design experiments over the instances of dataset C15 to obtain a set of values for MOSWACP parameters and optimal parameters for the proposed ERSA using response surface methodology (RSM). The RSM is a statistical and mathematical technique widely employed in experimental design and optimization processes \citep{khuri2010response}. It is a powerful tool for understanding complex relationships between input factors and the response of a system. The RSM simultaneously explores the optimal settings for multiple parameters, leading to improved efficiency and performance of a given process or system \citep{campatelli2014optimization}.
In this study, we input the parameters of ERSA and the instances of C15 to the RSM and determine optimal parameters according to the objective function values.

\begin{figure}[!ht]
    \centering
    \fbox{\includegraphics[scale=0.4]{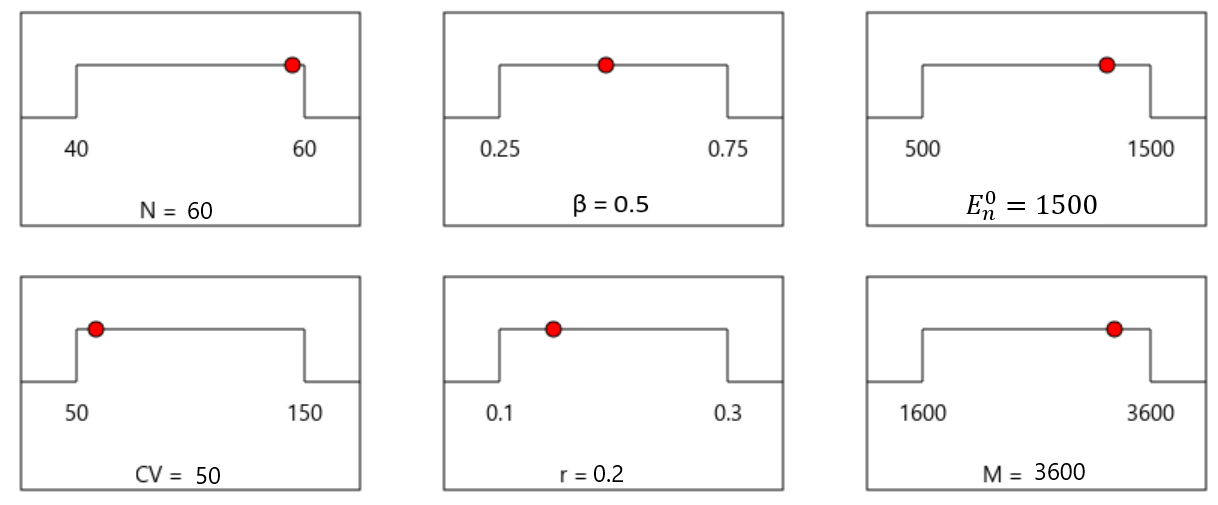}}
    \caption{Best values of ERSA parameters for MOSWACP based on the dataset C15}
    \label{fig:3}
\end{figure}

We illustrate the ranges and optimal values of the parameters obtained by ERSA in Fig. \ref{fig:3}. Since the output of RSM for the parameters, $N$, $E^0_n$, $CV$, and $M$ were float numbers, we rounded them to the nearest integer. Then, we found the best values for the parameters to minimize the objective function as $N=60$, $\beta=0.5$, $E^0_n=1500$, $CV=50$, $M=3600$, and $r=0.2$.
Also, Fig. \ref{fig:5} (see \ref{appendixA}) shows the behavior of the objective function value, $OF$ (represented on the y-axis), for different values of the parameters, $N$, $\beta$, $E^0_n$, $CV$, $r$, and $M$ (represented on the x-axis).
This figure shows that the objective function value negatively correlates with $M$ and $E^0_n$ but has no meaningful relationship with the remaining parameters.
The objective function value in these experiments was the average of the objective function value of the first $10$ instances of the C15 dataset returned by ERSA.

\subsection{Solutions to the adopted MIP model}

In this section, we solve MOSWACP at different sizes by the exact commercial solver, CPLEX, and illustrate computational results in Table \ref{table2}.

\begin{table}[ht!]
\centering
\scriptsize
\caption{The computational results on different datasets of the MOSWACP instances solved by CPLEX ($T$: execution time in seconds)}
\label{table2}
\begin{tabular}{l l l l p{1.5cm} p{1.5cm} p{1.5cm}}
\hline 
Dataset & Minimum T & Maximum T & Average T & \# executed instances & \# instances solved optimally & \% instances solved optimally\\
\hline
\textit{C15} & 5 & 2655 & 911 & 34 & 34 & 100\\
\textit{C21} & 15 & 3600 & 1100 & 84 & 77 & 92\\
\textit{J10} & 1 & 601 & 52 & 536 & 536 & 100\\
\textit{J12} & 1 & 3600 & 60 & 101 & 100 & 99\\
\textit{J14} & 4 & 3600 & 399 & 117 & 116 & 99 \\
\textit{J16} & 5 & 3600 & 498 & 45 & 41 & 91\\
\textit{J18} & 5 & 3600 & 501 & 93 & 68 & 73\\
\textit{J20} & 22 & 3600 & 1991 & 64 & 37 & 58\\
\textit{M2} & 3 & 3600 & 770 & 89 & 75 & 84\\
\textit{R3} & 5 & 3600 & 801 & 113 & 95 & 84 \\
\textbf{Average} & 6.60 & 3205.60 & 708.30 & — & — &88 \\
\hline 
\end{tabular}
\normalsize
\end{table}

Due to the execution time limit, we set it to one hour ($T=3600$ seconds), and left few instances from each dataset unsolved to optimality.
As shown in the table, we solved all instances of the datasets C15 and J10 in about $18$ and $1$ minutes, respectively.
Also, the MIP solver efficiently solves the instances of datasets C21, J12, J14, and J16 by obtaining optimal solutions in more than $90\%$ of the instances. However, solving the instances of dataset J20 is challenging for CPLEX since only less than $60\%$ of cases are solved optimally by CPLEX. On average, $88\%$ of all instances are solved optimally by CPLEX in about $12$ minutes.


\subsection{ERSA performance evaluation}

In this section, we solve MOSWACP instances of different sizes using our proposed ERSA, equipped with improvement rules. Also, we compare the performance of the proposed ERSA with SA, GA, and PSO, as summarized in Table \ref{table3}. We set the execution time limit to $1$ minute in these experiments and obtain the optimal value of parameters by RSM \citep{dean2017response}. We put the parameters of GA to $N=50, p_c=0.3, p_m=0.2, r_m=0.2, C_p=0.2$, where $N$ is the population size, $p_c$ is the crossover probability, $p_m$ is the mutation probability, $r_m$ is the mutation rate, and $C_p$ is the crossover point. Also, we found optimal values for the SA parameters as $T_{max}=1000, T_{min}=0.01, \alpha=0.98, N=20$, where $\alpha$ is the cooling rate, and $N$ is the number of iterations at each temperature. We also set the parameters of PSO to $n_s=20$, $M=100$, and $c_1=c_2=2$, where $n_s$ denotes the size of the swarm, $M$ is the maximum number of iterations. Also, $c_1$ and $c_2$ are the learning factors.

\begin{figure}[ht!]
\begin{tabular}{cc}
  \includegraphics[width=0.5\textwidth]{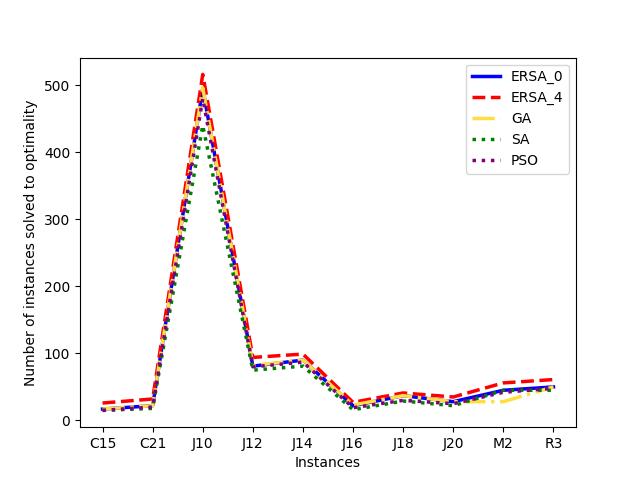} &
 \includegraphics[width=0.5\textwidth]{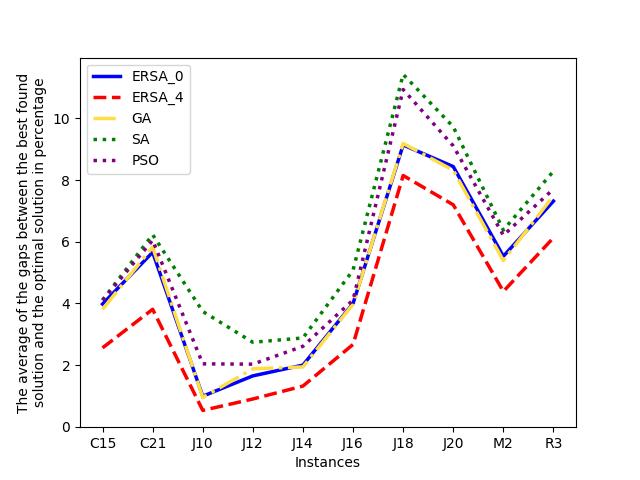} \\
 \footnotesize (a) The comparison of $N^*$ & \footnotesize (b) The comparison of $R^*$ \\[6pt]
\end{tabular}
\caption{The comparison of the enhanced ERSA with the existing solvers based on: (a) $N^*$, total instances solved to optimality, and (b) $R^*$, average of relative percentage between the best-found solutions}
\label{fig:comparison}
\end{figure}

In Table \ref{table3}, we illustrate $5$ different versions of ERSA equipped with a different set of improvement rules involved: (i) ERSA with no improvement rule, $ERSA_0$; (ii) ERSA with all improvement rules except $IR_3$, $ERSA_1$; (iii) ERSA with all improvement rules except $IR_2$, $ERSA_2$; (iv) ERSA with all improvement rules except $IR_1$, $ERSA_3$; and (v) ERSA with all improvement rules, $ERSA_4$. Also, we compare metaheuristic solvers according to three different criteria: (i) $N^*$, the number of instances solved to optimality; (ii) $R^*$, the average of relative percentage between the best-found solutions, $f$ (returned by metaheuristics), and the optimal solutions, $f^*$ (returned by MIP solver), where $R^*=\frac{f-f^*}{f^*}*100$; and (iii) $T^*$, the average execution time. As shown in Table \ref{table3}, our proposed ERSA with the proposed improvement rules performs better than existing metaheuristic solvers in terms of the number of instances solved to optimality and the average gap between the best-found and optimal solutions. Moreover, the designed improvement rules are superior to the solvers regarding solution quality and closeness to the optimal solution. We compare the proposed ERSA with the existing metaheuristic solvers in terms of the number of instances solved to optimality, $N^*$, as well as the average of relative percentage between the best-found solution and the optimal solution, $R^*$, in Fig. \ref{fig:comparison}. The result reveals that the proposed ERSA bundled with improvement rules, $ERSA_4$, has a better performance compared to the existing metaheuristic solvers.

\begin{table}[ht!]
\caption{The comparison of the proposed ERSA with the existing metaheuristics at various-sized instances of MOSWACP}
\label{table3}
\centering
\scriptsize
\scalebox{0.9}{
\begin{tabular}{c c c c c c c c c c c}
\hline  
Dataset & \# instances & Criteria & $ERSA_0$ & $ERSA_1$ & $ERSA_2$ & $ERSA_3$ & $ERSA_4$ & GA & SA & PSO \\
\hline
\textit{C15} & 34 & $N^*$ & 17 & 20 & 21 & 21 & 26 & 17 & 15 & 15 \\
 &  & $R^*$ & 3.98 & 3.74 & 3.72 & 3.24 & 2.56 &  3.81& 4.11 & 4.10 \\
 &  & $T^*$ & 31.09 & 35.55 & 36.13 & 39.15 & 43.19 & 39.01 & 45.11 & 36.09 \\
\textit{C21} & 84 & $N^*$ & 22 & 25 & 26 & 26 & 32 & 22 & 18 & 20 \\
 &  & $R^*$ & 5.66 & 5.34 & 5.33 & 5.21 & 3.81 &  5.83 & 6.23 & 6.05 \\
 &  & $T^*$ & 41.66 & 45.08 & 49.01 & 51.00 & 55.01 & 49.22 & 48.01 & 50.17 \\
 \textit{J10} & 536 & $N^*$ & 495 & 499 & 501 & 501 & 516 & 496 & 440 & 480 \\
 &  & $R^*$ & 0.99 & 0.83 & 0.81 & 0.80 & 0.53 &  0.95 & 3.74 & 2.04 \\
 &  & $T^*$ & 37.66 & 39.01 & 44.43 & 45.81 & 49.01 & 44.82 & 48.00 & 49.76 \\
  \textit{J12} & 101 & $N^*$ & 81 & 85 & 87 & 87 & 94 & 82 & 75 & 80 \\
 &  & $R^*$ & 1.65 & 1.50 & 1.32 & 1.12 & 0.90 &  1.88 & 2.74 & 2.03 \\
 &  & $T^*$ & 39.01 & 42.01 & 44.91 & 45.09 & 49.33 & 44.10 & 42.04 & 47.16 \\
 \textit{J14} & 117 & $N^*$ & 90 & 94 & 94 & 95 & 99 & 90 & 81 & 87 \\
 &  & $R^*$ & 2.00 & 1.78 & 1.50 & 1.45 & 1.32 &  1.94 & 2.88 & 2.61 \\
 &  & $T^*$ & 44.09 & 46.12 & 48.01 & 49.12 & 50.81 & 49.66 & 45.08 & 51.00 \\
 \textit{J16} & 45 & $N^*$ & 22 & 22 & 23 & 25 & 27 & 23 & 16 & 19 \\
 &  & $R^*$ & 4.00 & 3.41 & 3.30 & 3.01 & 2.67 &  3.98 & 5.08 & 4.11 \\
 &  & $T^*$ & 41.33 & 42.01 & 45.39 & 46.91 & 49.11 & 51.60 & 50.44 & 51.05 \\
 \textit{J18} & 93 & $N^*$ & 37 & 37 & 38 & 38 & 41 & 37 & 29 & 30 \\
 &  & $R^*$ & 9.14 & 8.21 & 8.20 & 8.18 & 8.15 &  9.19 & 11.43 & 10.94 \\
 &  & $T^*$ & 49.12 & 49.44 & 49.31 & 49.91 & 50.09 & 54.80 & 45.04 & 50.77 \\
\textit{J20} & 64 & $N^*$ & 28 & 28 & 29 & 33 & 35 & 28 & 22 & 25 \\
 &  & $R^*$ & 8.44 & 8.09 & 7.91 & 7.53 & 7.20 &  8.32 & 9.73 & 9.11 \\
 &  & $T^*$ & 42.22 & 43.01 & 45.90 & 47.38 & 49.00 & 49.32 & 44.01 & 48.71 \\
 \textit{M2} & 89 & $N^*$ & 45 & 48 & 50 & 54 & 56 & 28 & 45 & 42 \\
 &  & $R^*$ & 5.53 & 5.31 & 5.02 & 4.66 & 4.39 &  5.39 & 6.38 & 6.21 \\
 &  & $T^*$ & 42.06 & 45.71 & 46.03 & 48.91 & 50.33 & 56.81 & 42.01 & 48.97 \\
 \textit{R3} & 113 & $N^*$ & 50 & 53 & 55 & 58 & 61 & 50 & 45 & 48 \\
 &  & $R^*$ & 7.32 & 7.01 & 6.95 & 6.61 & 6.13 &  7.49 & 8.29 & 7.71 \\
 &  & $T^*$ & 46.11 & 48.64 & 48.93 & 48.99 & 52.94 & 55.99 & 49.22 & 53.01 \\
\hline
\end{tabular}}
\end{table}

\begin{figure}[ht!]
\begin{tabular}{ccc}
\includegraphics[width=0.3\textwidth]{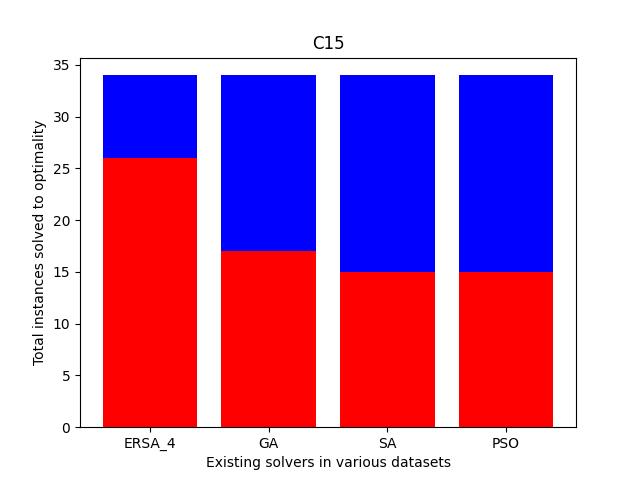} &
\includegraphics[width=0.3\textwidth]{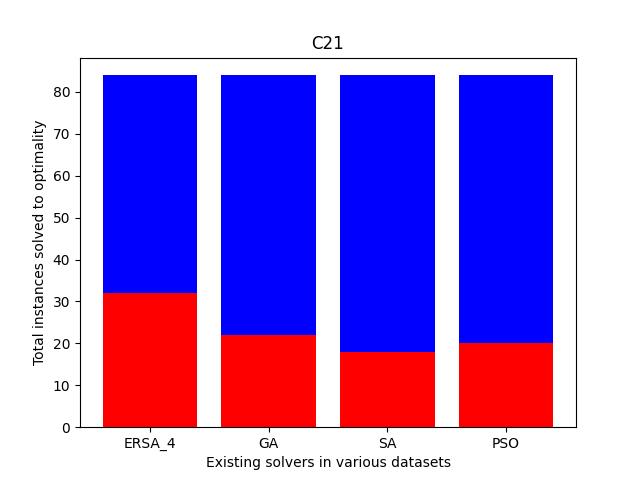} &
\includegraphics[width=0.3\textwidth]{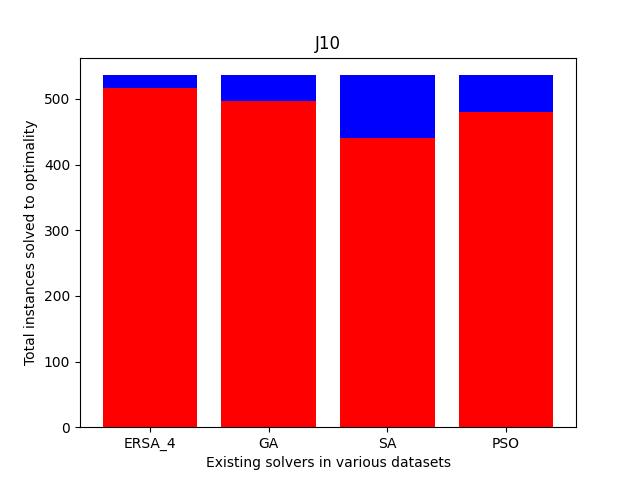} \\
\includegraphics[width=0.3\textwidth]{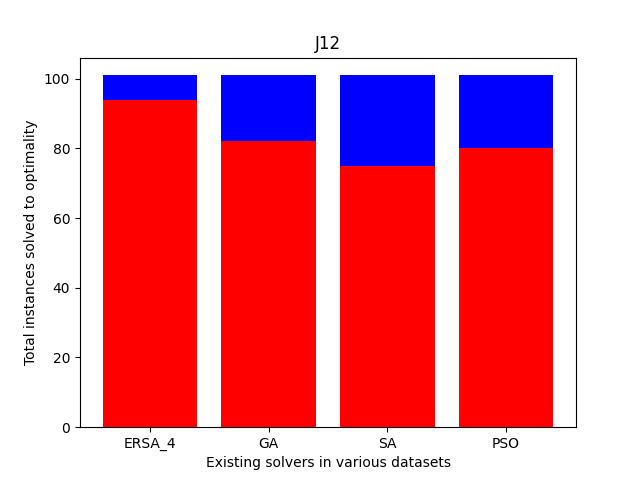} &
\includegraphics[width=0.3\textwidth]{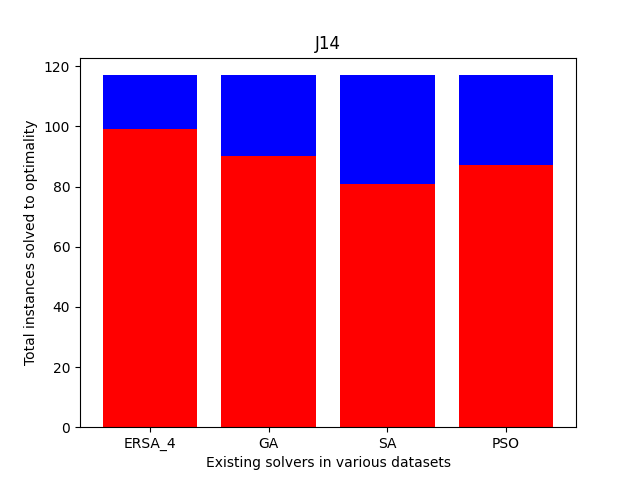} &
\includegraphics[width=0.3\textwidth]{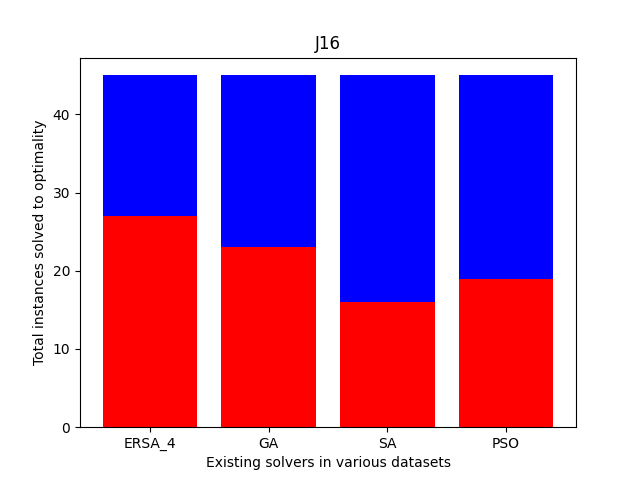} \\
\includegraphics[width=0.3\textwidth]{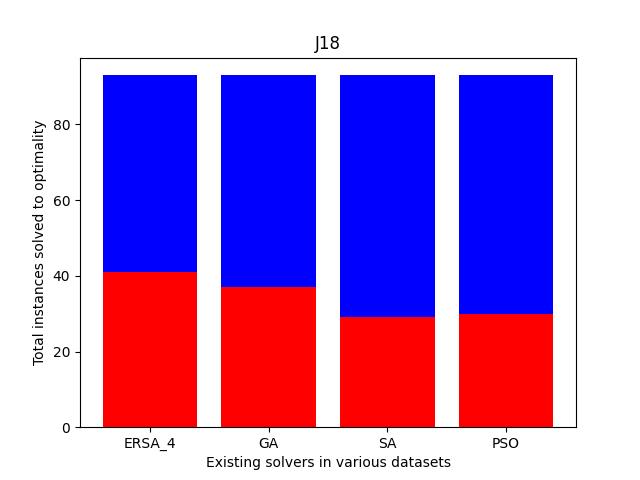} &
\includegraphics[width=0.3\textwidth]{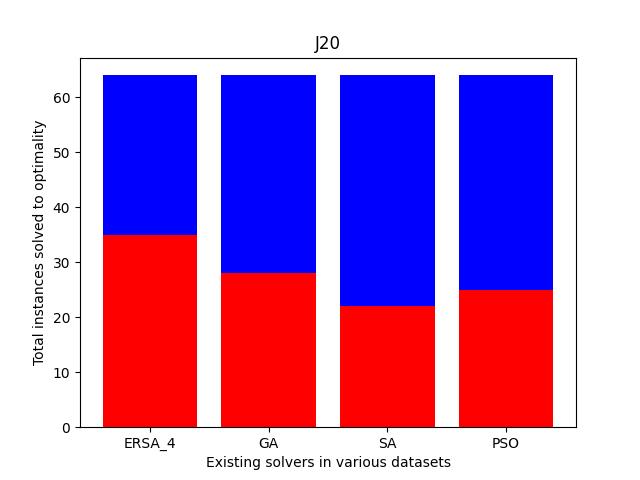} &
\includegraphics[width=0.3\textwidth]{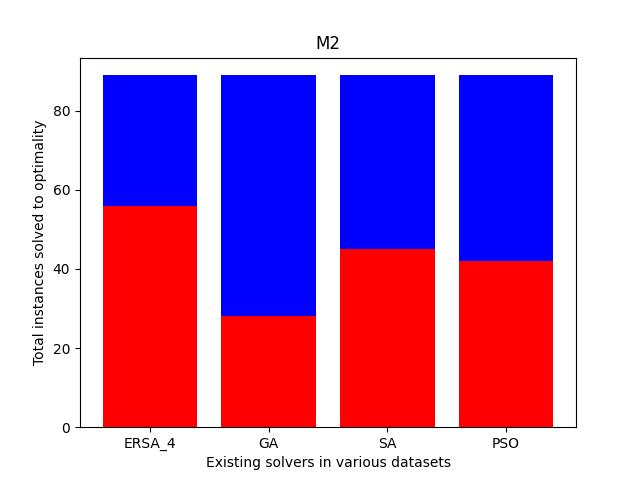} \\
\includegraphics[width=0.3\textwidth]{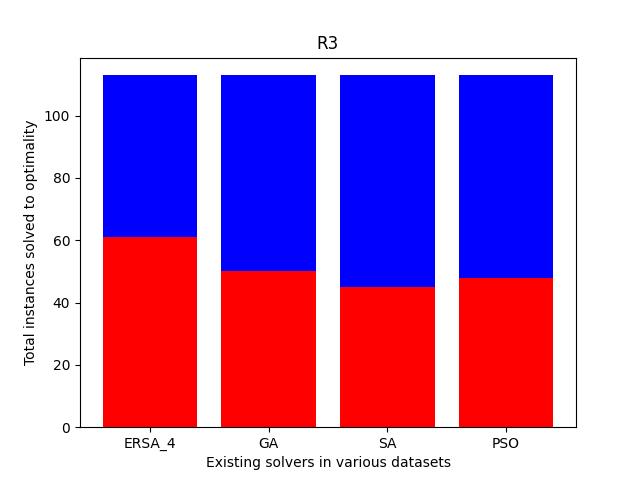} &
&
\\
\end{tabular}
\caption{The number of instances solved (blue) and not solved (red) to optimality in datasets}
\label{fig:red_blue}
\end{figure}

Also, we compare the fraction of the total instances solved to the optimality of each dataset. As shown in Fig.~\ref{fig:red_blue}, the cases solved optimally in each dataset are highlighted in blue, and those not solved optimally are shaded in red. The prominent blue areas in these figures demonstrate the enhanced ERSA, incorporating the proposed improvement rules and finding relatively better optimal solutions compared to the existing solvers.

\begin{figure}[ht!]
\begin{tabular}{cc}
\includegraphics[width=0.5\textwidth]{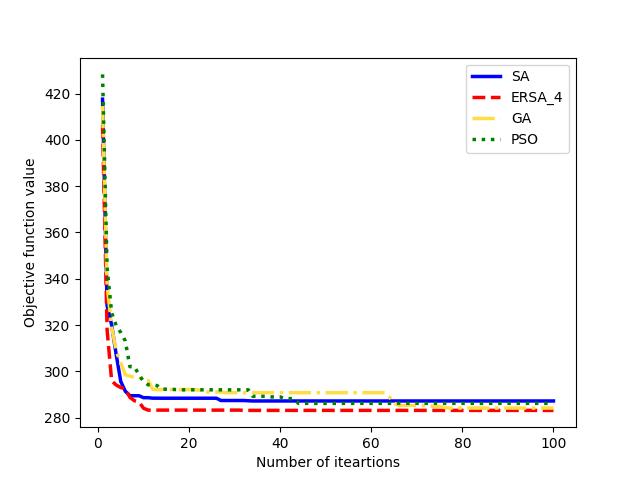} &
\includegraphics[width=0.5\textwidth]{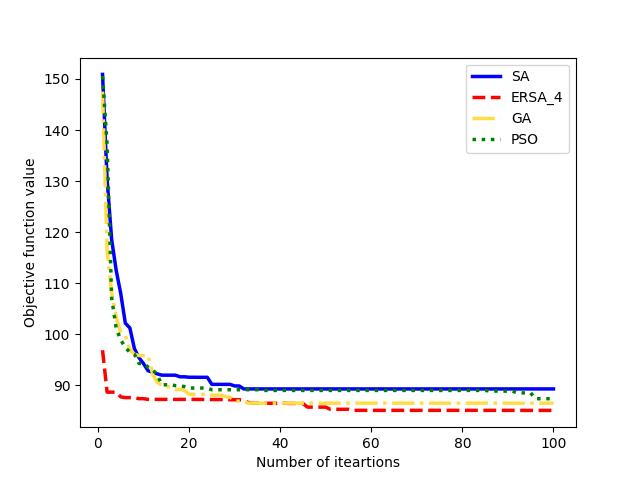} \\
\footnotesize (a) C21 & \footnotesize (b) J20 \\[6pt]
\includegraphics[width=0.5\textwidth]{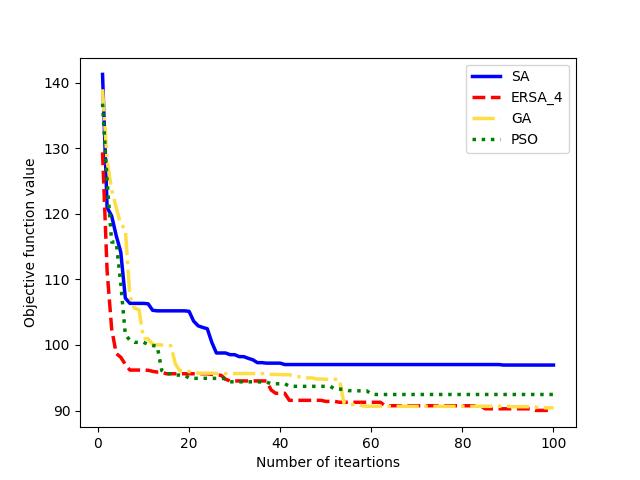} &
\includegraphics[width=0.5\textwidth]{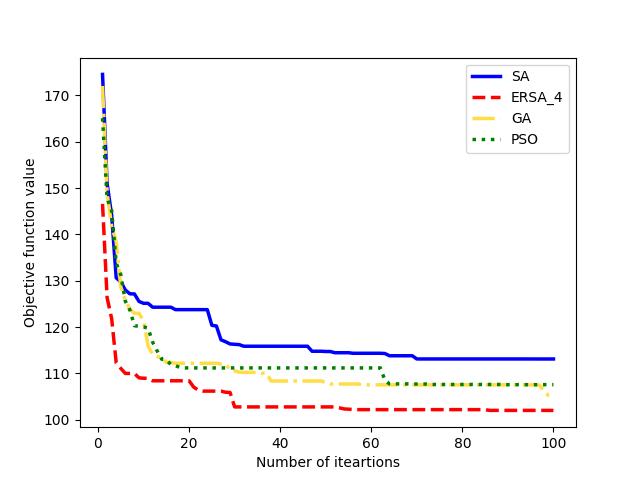}\\
\footnotesize (c) M2 & \footnotesize (d) R3
\end{tabular}
\caption{convergence history of $ERSA_4$, GA, PSO, and SA for instances of the datasets: (a) C21, (b) J20, (c) M2, and (d) R3}
\label{fig:convergence}
\end{figure}

Also, we compare the convergence of $ERSA_4$ with the other solvers in Fig. \ref{fig:convergence} for large instances of datasets C21, J20, M2, and R3, respectively.
Even though GA outperforms PSO and SA in finding relatively better solutions on average, our proposed $ERSA_4$ generates relatively higher-quality solutions. It searches for neighbor solutions more efficiently for large instances of MOSWACP.
These figures illustrate that $ERSA_4$ finds high-quality solutions in relatively earlier iterations and searches local solutions more efficiently than other methods.

Finally, we compare $ERSA_4$ with the CPLEX solver in the last instance of each dataset. As given in Table \ref{tableexact}, the proposed $ERSA_4$ reaches optimality in $2$ of datasets. Although the solutions by the CPLEX solver are acquired in one hour ($3600$ seconds), the proposed algorithm found solutions close to CPLEX solutions in approximately one minute ($60$ seconds) on average. Also, the gap, defined as the relative difference between the best-found solution by CPLEX and the proposed $ERSA_4$, depicts a faster performance of $ERSA_4$ for large instances of datasets compared to the CPLEX.

\begin{table}[ht!]
\centering
\footnotesize
\caption{Comparison of CPLEX and $ERSA_4$}\label{tableexact}

\scalebox{0.9}{
\begin{tabular}{cccc|ccc}
\hline
& \multicolumn{3}{c}{CPLEX}                              & \multicolumn{3}{|c}{$ERSA_4$}                           \\
Dataset  & Best-found solution & T & Optimal & Best-found solution & T & Gap (\%) \\
\hline
\textit{C15} & 103.45                       & 2655           & Yes              & 103.45                       & 63.64          & 0.00              \\
\textit{C21} & 280.98                       & 3600           & No               & 284.43                       & 61.44          & 1.23          \\
\textit{J10} & 43.38                        & 601            & Yes              & 43.38                        & 59.43          & 0.00              \\
\textit{J12} & 49.51                        & 3600           & NO               & 49.54                        & 60.94          & 0.06          \\
\textit{J14} & 57.74                        & 3600           & NO               & 58.43                        & 60.54          & 1.19          \\
\textit{J16} & 61.99                        & 3600           & NO               & 64.98                        & 62.81          & 4.82          \\
\textit{J18} & 71.04                        & 3600           & NO               & 75.23                        & 62.65          & 5.90          \\
\textit{J20} & 85.44                        & 3600           & NO               & 88.49                        & 63.28          & 3.57           \\
\textit{M2}  & 90.43                        & 3600           & NO               & 92.64                        & 63.88          & 2.44          \\
\textit{R3}  & 102.47                       & 3600           & NO               & 107.82                       & 62.91          & 5.22          \\
\hline
\textbf{Average}      & 94.643                       & 3205.6         &                  & 96.839                       & 62.15         & 2.44\\     \hline    
\end{tabular}}
\end{table}

The computational results demonstrate that our proposed Electron Radar Search Algorithm (ERSA) with all improvement rules outperforms the existing metaheuristics, achieving more instances solved to optimality and maintaining smaller gaps between the best-found solutions and optimal counterparts. Furthermore, the results highlight the effectiveness of the improvement rules compared to other solvers lacking enhancements, such as solution quality and proximity to the optimal solution.

\begin{figure}[ht!]
\centering
\fbox{\includegraphics[scale=0.7]{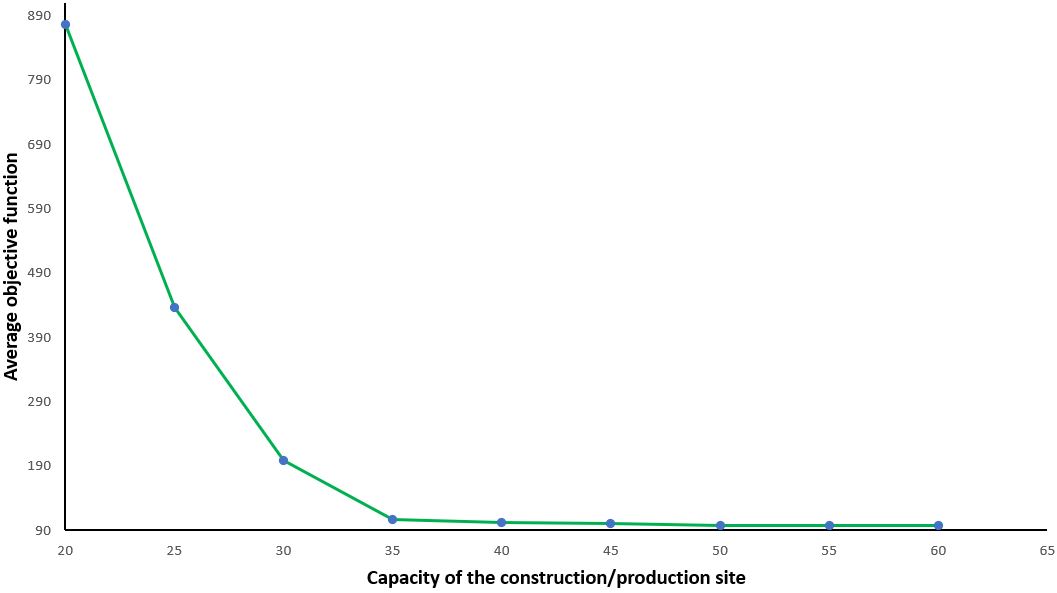}}
    \caption{Impact of project's site capacity ($Q$) on the usage cost of OSWs}
    \label{fig:sen-ana}
\end{figure}

Moreover, the parameter of the spatial capacity of the construction site (denoted by $Q$) is changed to observe its impact on the objective function. The ten instances of dataset C15 were selected, and the average objective function over these instances was reported. Fig. \ref{fig:sen-ana} shows the impact of changing the capacity of the construction/production site on the objective function (minimizing the usage cost of OSWs). As observed, the average cost is reducing (improving) by increasing the site's capacity ($Q$). The reason is that by providing more capacity for a construction site, the OSWs could be installed anytime. So, the associated activities could utilize these OSWs, leading to lower makespan and efficient utilization of OSWs. Also, in lower capacities, the marginal impact is higher; however, after some point, increasing the capacity does not change the usage cost of OSWs. Therefore, the project managers must know that the larger site's capacity does not constantly improve the costs, and the optimal value of the $Q$ could be found when the costs change drastically.  

\subsection{Case study}
This section presents the application of the OSWs and their practicability for an actual case study related to a trailer production project. In this case study, we aim to locate the OSWs within the production site (considering the occupied space and their lifetime) to make the project scheduling more efficient and to utilize the resources more effectively, compared to traditional policies. 
Detailed information on the company is confidential upon their request. Since 2000, the company has expanded production lines, parts making, and assembly to more than 30,000 square meters of land. Data from the case study may be available upon a reasonable request and owner's permission. The information on project activities, duration, and precedence relations related to this project of the trailer production line is shown in Table \ref{case-study-act}. Also, the parameters related to the OSWs and storage areas within the manufacturing site are shown in Table \ref{case-study-osw}.

\begin{table}[ht!]
\centering
\caption{Activities, their duration and precedence relations in the case study}\label{case-study-act}
\begin{tabular}{lll}
\hline
\multicolumn{1}{c}{\textbf{Activity}} & \multicolumn{1}{c}{\textbf{Duration (days)}} & \multicolumn{1}{c}{\textbf{Predecessors}} \\
\hline
A1 & 8 & - \\
A2                                    & 5                                            & A1                                        \\
A3                                    & 7                                            & A1                                        \\
A4                                    & 6                                            & A2, A3                                    \\
A5                                    & 10                                           & A4                                        \\
A6                                    & 5                                            & A4                                        \\
A7                                    & 12                                           & A5, A6                                    \\
A8                                    & 8                                            & A7\\
\hline
\end{tabular}
\end{table}

\begin{table}[ht!]
\centering
\footnotesize
\caption{OSWs, their duration and precedence relations in the case study}\label{case-study-osw}
\begin{adjustbox}{max width=\textwidth}
\begin{tabular}{llllll}
\hline
\multicolumn{1}{c}{\textbf{OSW}} & \multicolumn{1}{c}{\textbf{Cost per day (USD)}} & \multicolumn{1}{c}{\textbf{Installation Time (days)}} & \multicolumn{1}{c}{\textbf{Dismantling Time (days)}} & \multicolumn{1}{c}{\textbf{Space Occupied (units)}} & \multicolumn{1}{c}{\textbf{Max Lifetime (days)}} \\
\hline
W1                              & 300                                             & 2                                                     & 1                                                    & 1                                                   & 30                                               \\
W2                              & 400                                             & 3                                                     & 2                                                    & 1                                                   & 40                                               \\
W3                              & 250                                             & 1                                                     & 1                                                    & 2                                                   & 20                               \\
\hline
\end{tabular}
\end{adjustbox}
\end{table}

The project duration (deadline) is 50 days, with a total of eight activities (A1-A8) and three On-Site Workshops (W1-W3). The space limitation is two OSWs at a time. Also, OSWs must not exceed their maximum lifetime on-site. The project's purpose is to minimize total cost, which includes the cost of OSW availability while finding the installation and dismantling time of OSWs and the schedule of activities. Also, in the case study, one execution mode is considered for each activity. Also, the space occupied by each activity is assumed to be the same considering execution mode or OSW. The definition of each activity in the case study is stated as follows:

\begin{enumerate}[I]
    \item \textbf{A1: Chassis Preparation.} This is the initial phase where the trailer chassis (the base frame) is assembled. This involves cutting, welding, and assembling steel beams to create the skeleton of the trailer.
    \item \textbf{A2: Axle and Suspension Installation.} Once the chassis is ready, the axles and suspension system are installed. These components ensure that the trailer can bear weight and move smoothly. Its predecessor is A1 (Chassis Preparation).
    \item \textbf{A3: Wiring and Electrical Setup.} During this step, electrical wiring and systems are installed for brake lights, indicators, and internal electrical components. This includes routing cables along the chassis. Its predecessor is A1 (Chassis Preparation).
    \item \textbf{A4: Fabrication of Body Components.} The trailer body (e.g., side panels, floor, and roof) is fabricated and prepped. This involves metal cutting, shaping, and welding to create the trailer's exterior.
   Its predecessors are A2 (Axle and Suspension Installation), A3 (Wiring and Electrical Setup).
   \item \textbf{A5: Paint and Surface Treatment.} The entire trailer chassis and body are primed, painted, and treated with anti-corrosion materials to ensure durability and resistance to weather conditions.
   Its predecessor is A4 (Fabrication of Body Components).
   \item \textbf{A6: Interior Fitting and Insulation.} Insulation materials and internal components (like partitions, insulation layers, or specialized fittings) are installed inside the trailer if required (for refrigerated or specialized trailers).
   Its predecessor is A4 (Fabrication of Body Components).
   \item \textbf{A7: Assembly of External Components.} This stage involves assembling external features such as doors, loading ramps, handles, and other accessories. Quality checks are also conducted at this stage.
   Its predecessors are A5 (Paint and Surface Treatment), A6 (Interior Fitting and Insulation).
   \item \textbf{A8: Final Quality Control and Testing.} The trailer goes through final inspections, quality control checks, and testing (e.g., brake tests, durability tests, etc.) to ensure everything is functioning properly before delivery.
   Its predecessor is A7 (Assembly of External Components).
\end{enumerate}

After solving the case study with our proposed methodology, the optimal activities' start and finish times were obtained as given in Table \ref{case-study-opt-act}. Moreover, considering the spatial constraint (maximum of 2 OSWs installed at any time) and the lifetimes of the OSWs, the optimal schedule for OSW availability is shown in Table \ref{case-study-opt-osw}. According to the results, W1 is installed on Day 1 to handle initial activities like Chassis Preparation (A1) and Axle and Suspension Installation (A2). It also stays on-site until Day 26, covering the Interior Fitting (A6).
W2 is installed on Day 16 to handle body fabrication and painting tasks like Fabrication of Body Components (A4) and Paint and Surface Treatment (A5). It remains on-site for 27 days until Day 43, when the Assembly of External Components (A7) is completed. W3 is installed on Day 32 to handle the final stages of assembly and Final Quality Control (A8), being dismantled after Day 50. Its availability is shorter, as it serves only during the final phase. The total usage cost for the OSWs is computed based on their availability times and per-day costs:

\begin{equation*}
\text{W1 Cost: } 25 \text{ days} \times 300 \text{ USD/day} = 7,500 \text{ USD} 
\end{equation*}

\begin{equation*}
\text{W2 Cost: } 27 \text{ days} \times 400 \text{ SD/day} = 10,800 \text{ USD}  
\end{equation*}

\begin{equation*}
\text{W3 Cost: } 18 \text{ days} \times 250 \text{ USD/day } = 4,500 \text{ USD} 
\end{equation*}

Thus, the total OSW usage cost (objective function value) is:

\begin{equation*}
\text{Total Cost}=7,500
+10,800+4,500=22,800 \text{ USD} 
\end{equation*}

Activities' optimal start and finish times are provided, ensuring all precedence constraints are respected. The OSWs are optimally scheduled with installation and dismantling times to minimize costs while adhering to spatial constraints. The total usage cost of the OSWs (objective function) is 22,800 USD. This setup ensures that the resources are optimally allocated, minimizing the costs associated with the OSWs.

\begin{table}[ht!]
\centering
\footnotesize
\caption{Optimal schedule of activities}\label{case-study-opt-act}
\begin{adjustbox}{max width=\textwidth}
\begin{tabular}{llll}
\hline
\multicolumn{1}{c}{\textbf{Activity}} & \multicolumn{1}{c}{\textbf{Duration (days)}} & \multicolumn{1}{c}{\textbf{Start Time (day)}} & \multicolumn{1}{c}{\textbf{Finish Time (day)}} \\
\hline
A1                                    & 8                                            & Day 1                                         & Day 8                                          \\
A2                                    & 5                                            & Day 9                                         & Day 13                                         \\
A3                                    & 7                                            & Day 9                                         & Day 15                                         \\
A4                                    & 6                                            & Day 16                                        & Day 21                                         \\
A5                                    & 10                                           & Day 22                                        & Day 31                                         \\
A6                                    & 5                                            & Day 22                                        & Day 26                                         \\
A7                                    & 12                                           & Day 32                                        & Day 43                                         \\
A8                                    & 8                                            & Day 44                                        & Day 50        \\
\hline
\end{tabular}
\end{adjustbox}
\end{table}

\begin{table}[ht!]
\centering
\caption{Optimal installation and dismantling times of OSWs in the case study}\label{case-study-opt-osw}
\begin{adjustbox}{max width=\textwidth}
\begin{tabular}{llllll}
\hline
\multicolumn{1}{c}{\textbf{OSW}} & \multicolumn{1}{c}{\textbf{Installation Time (day)}} & \multicolumn{1}{c}{\textbf{Dismantling Time (day)}} & \multicolumn{1}{c}{\textbf{Availability Duration (days)}} & \multicolumn{1}{c}{\textbf{Space Occupied (units)}} & \multicolumn{1}{c}{\textbf{Total Usage Cost (USD)}} \\
\hline
W1& Day 1& Day 26   & 25& 1 & 300 * 25 = 7,500\\
W2 & Day 16& Day 43     & 27& 1 & 400 * 27 = 10,800                                   \\
W3                               & Day 32                                               & Day 50                                              & 18                                                        & 2                                                   & 250 * 18 = 4,500                                   \\
\hline
\end{tabular}
\end{adjustbox}
\end{table}

In the following, we compare the total cost returned by our methodology with the total cost obtained by the traditional policies (based on past experiences). The activity's start and finish times and OSW's dismantling and installation time returned by the conventional policy are in Tables \ref{case-study-policy-act}-\ref{case-study-policy-osw}, respectively. Based on these tables, the total usage cost of OSWs when using traditional policy equals $9,000 + 16,800 + 8,750 = 34,550$ USD. Compared with the optimal solution (22,800 USD), the savings with the Optimal Solution equal $34,550 - 22,800 = 11,750$ USD. So, the optimal solution is 33.99\% cheaper than the traditional approach. The reason is that in the optimal solution, the installation and dismantling of OSWs are timed more efficiently, reducing idle times. On the other hand, in the traditional solution, the OSWs are installed too early and kept for too long, which leads to unnecessarily higher usage costs.

\begin{table}[ht!]
\centering
\footnotesize
\caption{The start and finish times of activities in the case study using traditional policy}\label{case-study-policy-act}
\begin{adjustbox}{max width=\textwidth}
\begin{tabular}{llll}
\hline
\multicolumn{1}{c}{\textbf{Activity}} & \multicolumn{1}{c}{\textbf{Duration (days)}} & \multicolumn{1}{c}{\textbf{Start Time (day)}} & \multicolumn{1}{c}{\textbf{Finish Time (day)}} \\
\hline
A1                                    & 8                                            & Day 1                                         & Day 8                                          \\
A2                                    & 5                                            & Day 9                                         & Day 13                                         \\
A3                                    & 7                                            & Day 9                                         & Day 15                                         \\
A4                                    & 6                                            & Day 16                                        & Day 21                                         \\
A5                                    & 10                                           & Day 22                                        & Day 31                                         \\
A6                                    & 5                                            & Day 22                                        & Day 26                                         \\
A7                                    & 12                                           & Day 32                                        & Day 43                                         \\
A8                                    & 8                                            & Day 44                                        & Day 50        \\
\hline
\end{tabular}
\end{adjustbox}
\end{table}

\begin{table}[ht!]
\centering
\caption{The installation and dismantling times of OSWs in the case study using traditional policy}\label{case-study-policy-osw}
\begin{adjustbox}{max width=\textwidth}
\begin{tabular}{llllll}
\hline
\multicolumn{1}{c}{\textbf{OSW}} & \multicolumn{1}{c}{\textbf{Installation Time (day)}} & \multicolumn{1}{c}{\textbf{Dismantling Time (day)}} & \multicolumn{1}{c}{\textbf{Availability Duration (days)}} & \multicolumn{1}{c}{\textbf{Space Occupied (units)}} & \multicolumn{1}{c}{\textbf{Total Usage Cost (USD)}} \\
\hline
W1                               & Day 1                                                & Day 31                                              & 30                                                        & 1                                                   & 300 * 30 = 9,000                                    \\
W2                               & Day 1                                                & Day 43                                              & 42                                                        & 1                                                   & 400 * 42 = 16,800                                   \\
W3                               & Day 16                                               & Day 50                                              & 35                                                        & 2                                                   & 250 * 35 = 8,750                                   \\
\hline
\end{tabular}
\end{adjustbox}
\end{table}

Furthermore, the optimal solution respects the space constraint by installing only two OSWs simultaneously, ensuring minimal overlap and avoiding resource congestion. The solution using traditional policy doesn't fully account for space optimization. Although it respects the constraint, the inefficiency in scheduling leads to wasted space, as some OSWs remain idle for extended periods. Also, the optimal solution leverages more precise timing of the availability of OSWs, coordinating their use exactly when needed and dismantling them when their role is complete. In contrast, the traditional solution installs OSWs too early and dismantles them too late, resulting in longer usage durations and higher costs.

In addition, in the optimal solution, OSWs are scheduled to avoid idle periods, ensuring they are only used when necessary for the associated activities. The traditional solution, however, leads to extended periods where OSWs are idle, occupying space but not contributing to the project, leading to unnecessary costs. Regarding cost efficiency, the optimal solution minimizes the number of days that OSWs are in use, resulting in a significant cost reduction. In the traditional solution, the usage periods are excessive, inflating the cost by nearly 12,000 USD.

\section{Discussion and Conclusion} \label{sec:conclusions}
When we compare ERSA with improvement rules ($ERSA_4$) to the existing metaheuristic solvers using two metrics, namely $N^*$ (total instances solved to optimality) and $R^*$ (the average difference (relative percentage) between the best-found solution and the optimal solution), the $ERSA_4$ exhibits superior efficiency. One notable strength of $ERSA_4$ is the rapid convergence to the best-found solution in early iterations, which sets it apart from alternative methods. The $ERSA$ excels at generating and exploring neighboring solutions, ultimately leading to the exploration of high-quality solutions. With large instances of MOSWACP, the Genetic Algorithm (GA) outperforms Particle Swarm Optimization (PSO) and Simulated Annealing (SA). GA consistently achieves solutions with lower objective function values, indicating its prowess in finding superior solutions.

Moreover, as we compared the solutions generated by the CPLEX solver and $ERSA_4$, our proposed algorithm achieves optimal solutions in $2$ datasets and, on average, reaches the solutions close to CPLEX solutions in $8$ other instances by approximately $2.34\%$. This paper adopted the project resource cost availability problem with OSWs and multi-mode activities, a mixed-integer programming model \citep{moradi2023site}. We developed a model-specific Electron Radar Search Algorithm (ERSA) with purposeful improvement rules to solve large-scale MOSWACP. The proposed method consistently outperforms other metaheuristic solvers by achieving relatively high-quality solutions and maintaining smaller gaps between the best-found and optimal solutions, specifically in early iterations in relatively short times. An actual case study related to the trailer production site was solved by the proposed methodology. The optimal solution was significantly better due to its focus on timing and efficient resource management. It resulted in 11,750 USD in savings (a 33.99\% reduction in costs) compared to traditional policies. This highlights the importance of optimization in resource allocation and scheduling, which is achieved effectively using the proposed model.

For future direction, we recommend focusing on the location of workshops since their neighborhood and other resources might matter in real-world practices. Also, MOSWACP, with material procurement limitations, is a potential opportunity to develop. We also suggest extending existing metaheuristics to solve MOSWACP and comparing them with benchmark solvers. Also, uncertain and stochastic parameters related to resource availability could be addressed by agent-based simulation models \citep{aftabi2025sd}. Such exploration helps the construction industry deploy superior methods to schedule projects optimally and use large-scale solvers for project initiatives.



\section*{Consent for publication}
All authors consent to the publication. All authors read and approved the final manuscript.

\section*{Funding} The authors received no financial support for this paper's research, authorship, and publication.

\section*{Data Availability Statement}Data, models, and codes are available upon request.



\section*{Conflict of interest} The authors have no conflicts of interest to disclose.

\section*{Acknowledgements}
Not applicable.

\bibliographystyle{elsarticle-harv} 
\bibliography{references}

\appendix

\setcounter{figure}{0}    

\section{Objective function behavior with respect to ERSA parameters} \label{appendixA}

\begin{figure}[ht!]
    \centering
    \includegraphics[scale=0.7]{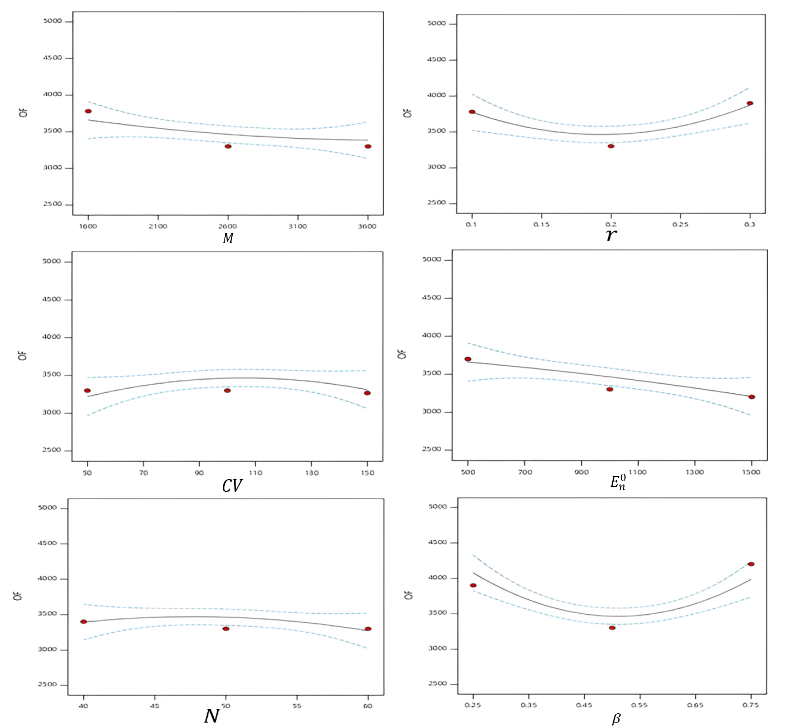}
    \caption{The behavior of the objective function ($OF$) with respect to the parameters of ERSA}
    \label{fig:5}
\end{figure}

Fig. \ref{fig:5} illustrates how the objective function value (OF) changes at different values of the parameters, $N, \beta, E^0_n, CV, r$, and $M$, on the vertical and horizontal axis, respectively.
The objective function value in these experiments is an average derived from ERSA's first $10$ instances of the dataset, C15.

\end{document}